\documentclass[lettersize,journal]{IEEEtran}
\pdfoutput=1

\usepackage[utf8]{inputenc}
\usepackage[T1]{fontenc}
\usepackage{graphicx}
\usepackage{grffile}
\usepackage{longtable}
\usepackage{wrapfig}
\usepackage{rotating}
\usepackage[normalem]{ulem}
\usepackage{amsmath}
\usepackage{textcomp}
\usepackage{amssymb}
\usepackage{capt-of}
\usepackage{hyperref}
\usepackage{amsthm}
\usepackage{mathtools}
\usepackage{centernot}
\usepackage{textgreek}
\usepackage{subcaption}
\usepackage{cite}

\hypersetup{colorlinks=true}
\newcommand{\real}{\mathbb{R}}
\newcommand\df{\stackrel{\mathclap{\tiny\mbox{df.}}}{=}}
\newcommand{\sgens}[1]{\mathbb{R}^{#1}}
\newcommand{\sgen}[2]{\mathbb{R}^{#1\times #2}}

\newcommand{\op}[1]{\operatorname{#1}}

\newcommand{\Int}[1]{\operatorname{int}({#1})}
\newcommand{\Intr}{\operatorname{int}({\sgens{k}_+})}
\newcommand{\Intrs}{\operatorname{int}({\mathbb{R}_+})}

\newcommand{\restr}[1]{_{\mkern 1mu \vrule height 2ex\mkern2mu{#1}}}
\newtheorem{assumption}{Assumption}
\newtheorem{example}{Example}
\newtheorem{proposition}{Proposition}
\newtheorem{fact}{Fact}
\newtheorem{remark}{Remark}

\newtheorem{definition}{Definition}

\newtheorem{lemma}{Lemma}

\newtheorem{corollary}{Corollary}

\begin{document}
\title{The fixed point iteration of positive concave mappings converges geometrically if a fixed point exists: implications to wireless systems}
\author{Tomasz Piotrowski$^{\dagger}$ and Renato L. G. Cavalcante$^{\dagger,\star}$\\
  $^1$ Faculty of Physics, Astronomy and Informatics,\\
  Nicolaus Copernicus University,
  Grudziądzka 5/7, 87-100 Toruń, Poland\\
  {$^2$ Fraunhofer Heinrich Hertz Institute \\
  Einsteinufer 37,  10587 Berlin, Germany\\}
}

\maketitle

\begin{abstract}
	We prove that the fixed point iteration of arbitrary positive concave mappings with nonempty fixed point set converges geometrically for any starting point. We also show that positivity is crucial for this result to hold, and the concept of (nonlinear) spectral radius of asymptotic mappings provides us with information about the convergence factor. As a practical implication of the results shown here, we rigorously explain why some power control and load estimation algorithms in wireless networks, which are particular instances of a fixed point iteration, have shown geometric convergence in simulations. These algorithms have been typically derived by considering fixed point iterations of the general class of standard interference mappings, so the possibility of sublinear convergence rate could not be ruled out in previous studies, except in special cases that are often more restrictive than those considered here. 
\end{abstract}


$\dagger$ Both authors contributed equally.

$\star$ Corresponding author.

\IEEEpeerreviewmaketitle

\section{Introduction}

Many problems in economy \cite{woodford2020}, wireless communications \cite{yates95,Cavalcante2016,Cavalcante2019_IEEE_TSP,martin11,slawomir09,you2020note,shindoh2020,cavalcante2016low,ho2014data}, and machine learning \cite{piotrowski2021fixed}, to cite a few fields, use the so-called \emph{standard interference mappings} \cite{yates95} in their formulations, and we note that these mappings are also known as (positive) order-preserving (or monotone) strictly subhomogeneous mappings in the mathematical literature \cite{Lemmens2012,Oshime1992,krause2015positive}. If we restrict the attention to the wireless domain, problems such as load estimation \cite{Cavalcante2016,Cavalcante2019_IEEE_TSP,Majewski2010,siomina12,feh2013,ho2014data}, power control \cite{martin11,slawomir09,you2020note,cavalcante2016low}, and signal-to-interference-noise ratio (SINR) feasibility studies \cite{shindoh2020,shindoh2019} often involve standard interference mappings. From a mathematical perspective, these problems have in common  that the objective is to answer whether a standard interference mapping has a fixed point, and, if so, whether the fixed point, which is unique if it exists \cite{yates95}, can be obtained with simple algorithms.  

Regarding the first question, the existence of a fixed point, recently the study in \cite{Cavalcante2019_IEEE_TSP} has used the concept of spectral radius of asymptotic mappings to obtain a necessary and sufficient condition for an arbitrary (possibly nonlinear) standard interference mapping to have a fixed point (see Fact \ref{sr1}). This condition is often easy to verify in practice, and it generalizes and unifies many mathematical tools for feasibility studies in the wireless literature. If a standard interference mapping has a fixed point, then its usage with the fixed point iteration produces a sequence that converges to the fixed point of the mapping \cite{yates95}, thus we have an answer to the second question posed above; namely, how to obtain the fixed point. 

Given the simplicity and widespread use of the fixed point iteration of standard interference mappings in the wireless domain, a great deal of effort has been devoted to understanding and improving its convergence rate  \cite{Cavalcante2016,cavalcante2018spectral,fey2012,huang1998rate}\cite[Ch.~5.3]{martin11}. From a theoretical perspective, an important but somewhat disappointing result in this direction has been established in \cite[Example~2]{fey2012}, which shows that the fixed point iteration of standard interference mappings can converge sublinearly. However, this bad performance does not usually manifest in practice. A possible explanation is that very frequently the standard interference mappings have additional structure that is not exploited to prove convergence of the fixed point iteration. In particular, positive concave mappings are common in applications \cite{yates95,huang1998rate,boche2008superlinearly,martin11,slawomir09,Cavalcante2016,Cavalcante2019_IEEE_TSP,cavalcante2016low,ho2014data,hanly1995algorithm}, and we recall that these mappings are a subclass of standard interference mappings \cite[Proposition~1]{Cavalcante2016}. Nevertheless, to the best our knowledge, to date it is unknown whether this structure is sufficient to guarantee geometric convergence in normed vector spaces except in special cases. There are recent results proving geometric convergence of fixed point iterations of a very general class of mappings \cite{luke2020necessary,Lauster2021}, but whether these results give rise to easily identifiable subclasses of standard interference mappings that guarantee geometric convergence of the fixed point iterations (and obtaining easily computable bounds on the convergence rate) remains largely unknown. 

Geometric convergence in normed vector spaces has been established if the fixed point iteration uses a concave mapping constructed by taking the coordinate-wise minimum of \emph{finitely} many positive affine mappings \cite{huang1998rate}, provided that a fixed point exists. Nevertheless, extending the analysis in that study to arbitrary positive concave mappings does not seem trivial, and we emphasize that concave mappings without the particular structure considered in \cite{huang1998rate} are very common in the wireless domain \cite{Cavalcante2016,Cavalcante2019_IEEE_TSP,cavalcante2016low,Majewski2010,siomina12,feh2013,boche2008,boche2008superlinearly}. In \cite[Example~2]{boche2008superlinearly}, the authors show that the fixed point iterations of positive affine mappings (a particular case of positive concave mappings) cannot converge super-linearly, and in \cite[Sect. VI]{boche2008superlinearly} they illustrate the convergence speed of a nonlinear positive concave mapping with figures. However, that discussion is restricted to a very specific nonlinear mapping related to a two-user power control problem with beamforming. Even in that specific setting, no formal proof of geometric convergence is provided. 

If a positive concave mapping has a fixed point, the study in \cite{cavalcante2018spectral} has shown that, with a mild assumption, the convergence of the fixed point iteration cannot be better than geometric, and, for specific starting points, the convergence factor is bounded by the spectral radius of the asymptotic mapping associated with the concave mapping. Nevertheless, those results do not rule out the possibility of sublinear convergence. In this study, we give a complete answer to this question. In more detail, our main contributions can be summarized as follows:

\begin{enumerate}

\item We prove that the fixed point iteration of an arbitrary (continuous)\footnote{The assumption of continuity can be dropped in many of our results, but we only consider mappings that are everywhere continuous on their domains to avoid technical digressions. We recall that positive concave mappings can be discontinuous on the boundary of their domains. See \cite{Cavalcante2016} for a simple example.} positive concave mapping with an arbitrary starting point is guaranteed to converge geometrically in normed vector spaces if the mapping has a fixed point (Proposition \ref{p1}). We also construct an example illustrating that this property does not necessarily hold if we drop the assumption of positivity (Example~\ref{renato_ex}).

\item We establish useful connections between the convergence factor of the fixed point iteration in Thompson's metric space and normed vector spaces (Proposition \ref{conv}).

\item We show that the concept of (nonlinear) spectral radius of the asymptotic mapping associated with a continuous positive concave mapping provides us with useful information about the convergence factor for any starting point of the fixed point iteration (Proposition~\ref{prop.linearc}).
  
\item We illustrate the main implications of the results derived here in power control and load estimation problems in wireless networks (Sect.~\ref{applications} and Sect.~\ref{simulations})
\end{enumerate}

To keep this study self-contained, we list in Appendix \ref{kru} definitions and facts that are required in proofs but not for understanding the main results of this study. Technical prerequisites that are crucial for setting the stage for the main contributions are discussed in the next section. In a first reading, we invite readers interested in problems in wireless networks to skip directly to Sect.~\ref{applications}  to understand the practical implications of the main technical contributions in Sects.~\ref{preliminaries} and \ref{main}. 

\section{Preliminaries}
\label{preliminaries}
Before introducing definitions and technical results required for the main proofs, we need the following necessary notation on the nonnegative cone:
\begin{equation*}
\sgens{k}_+\df\{(x_1,\dots,x_k)\in\sgens{k}:\ x_i\geq 0\textrm{ for }1\leq i\leq k\}.
\end{equation*}
By $\op{int}(\sgens{k}_+)$ we denote the set of (strictly) positive vectors, i.e., $x=(x_1,\dots,x_k)\in\op{int}(\sgens{k}_+)\iff x_i>0$ for $1\leq i\leq k.$ Let now $x,y\in\sgens{k}_+$; then, $x\leq y$ denotes the partial ordering induced by the nonnegative cone, i.e., $x \le y \iff y-x\in\sgens{k}_+$. Similarly, $x<y$ means that $y-x\in\sgens{k}_{+}$ and $x\neq y$, while $x\ll y$ means that $y-x\in\op{int}(\sgens{k}_+)$.

Let $x$ be an element of a topological space. We call $U$ a neighborhood of $x$ if there exists an open set $V$ such that $x\in V\subseteq U.$ We denote by $x[i]\in\mathbb{R}$ for $i\in\{1,2,\dots,k\}$ the $i$-th coefficient of a vector $x\in\sgens{k}$, and by $[X]_{i,j}\in\mathbb{R}$ the element of a matrix $X\in\sgen{m}{n}$ at the $i$-th row and $j$-th column with $i\in\{1,2,\dots,m\}$ and $j\in\{1,2,\dots,n\}.$

Given a mapping $f:X\to Y$ for two sets $X$ and $Y$ satisfying $Y\subseteq X\subseteq \sgens{k}_+$, we denote by
\begin{equation*}
\mathrm{Fix}(f)\df\{x\in X~|~f(x)=x\}
\end{equation*}
the set of fixed points of $f$. Moreover, for a metric space $(X,d_1)$, a sequence $(x_n)_{n\in\mathbb{N}}$ in $X$ converges to $\widetilde{x}\in X$ if $\forall\epsilon>0$ $\exists n_0\in\mathbb{N}$ such that $\forall n>n_0$ one has that $d_1(x_n,\widetilde{x})<\epsilon.$ We also note that a norm $\|\cdot\|$ of a normed vector space $(\sgens{k}_+,\|\cdot\|)$ induces a metric on $X$ by $d_2(x,y)\df\|x-y\|$ for $x,y\in X.$ In such a case, we say that $d_2$ is induced by a norm $\|\cdot\|.$

	Let $f\colon\sgens{k}_+\to\sgens{k}_+$, and construct a sequence $(x_n)_{n\in\mathbb{N}}$ via
	\begin{equation} \label{fpseq}
		\forall n\in\mathbb{N}\quad x_{n+1}\df f(x_n)\textrm{ with }x_1\in\sgens{k}_+.
	\end{equation}
	We call the recursion in (\ref{fpseq}) {\it a fixed point iteration} of $f$ because (\ref{fpseq}) is a classical algorithm to compute a point in $\mathrm{Fix}(f)$ provided that $\mathrm{Fix}(f)\neq\emptyset$. 
	
	Under consideration in this study is the problem of establishing geometric convergence, in the sense defined below, of a sequence constructed with the iteration in (\ref{fpseq}) with a particular class of mappings described later in Definition~\ref{PC}.

\begin{definition} \label{contr}
	Let $(x_n)_{n\in\mathbb{N}}$ be a sequence in $\sgens{k}_+.$ We say that $(x_n)_{n\in\mathbb{N}}$ converges geometrically to $\widetilde{x}\in\sgens{k}_+$ with a factor $c\in[0,1)$ and a constant $\gamma>0$ if
	\begin{equation} \label{geom}
		\exists c\in[0,1)\quad\exists \gamma>0\quad\forall n\in\mathbb{N}\quad \|x_{n+1}-\widetilde{x}\|\leq\gamma c^n,
	\end{equation}
	where $\|\cdot\|$ is an arbitrary norm. Similarly, we say that $(x_n)_{n\in\mathbb{N}}$ converges linearly to $\widetilde{x}\in\sgens{k}_+$ with a factor $c\in[0,1)$ if
	\begin{equation} \label{linear}
		\exists c\in[0,1)\quad\forall n\in\mathbb{N}\quad \|x_{n+1}-\widetilde{x}\|\leq c\,\|x_n-\widetilde{x}\|,
	\end{equation}
	where $\|\cdot\|$ is an arbitrary norm. We verify that, if $(x_n)_{n\in\mathbb{N}}$ converges linearly to $\widetilde{x}\in\sgens{k}_+$ with a factor $c\in[0,1)$, then it converges geometrically to $\widetilde{x}\in\sgens{k}_+$ with the same factor $c$ and a constant $\gamma=\|x_1-\widetilde{x}\|.$\footnote{Geometric convergence is also known in literature as R-linear convergence, whereas linear convergence is also known as Q-linear convergence \cite[Chapter 9]{Ortega1970}.} We also note that the above definitions of geometric and linear convergence in normed spaces can be naturally generalized to metric spaces, with the metric induced by the norm in (\ref{geom}) and (\ref{linear}) replaced with a metric of a  metric space.
\end{definition}

Given an arbitrary mapping $f\colon\sgens{k}_+\to\sgens{k}_+$, in general there are no guarantees that the fixed point iteration converges, even if $f$ has a fixed point. As a result, we need to impose additional structure on $f$, and the structure should keep the class of mappings general enough to be able to address as many real-world problems as possible. Of particular interest in the wireless literature are the standard interference mappings, which are known to produce fixed point iterations that converge if a fixed point exists \cite{yates95}. These mappings are defined below.

\begin{definition}
	\label{def.si}
A mapping $f\colon\sgens{k}_+\to\sgens{k}_+$ is said to be a standard interference mapping (SI mapping) if it is 
 \begin{enumerate}
 	\item \emph{monotonic}
 	\begin{equation}
 	\forall x,y\in\sgens{k}_+\quad x\leq y\implies f(x)\leq f(y),\textrm{ and }
 	\end{equation}
        \item \emph{scalable}
 	\begin{equation} \label{scala}
 	\forall x\in\sgens{k}_+\quad \forall \lambda> 1\quad f(\lambda x)\ll\lambda f(x).
 	\end{equation} 	
 \end{enumerate}
\end{definition}

\begin{remark}
	Functions satisfying the monotonicity property in Definition~\ref{def.si} receive different names depending on the author. In particular, they are also known as order-preserving. Hereafter, we simply call them monotone or monotonic functions, as also done in some studies on mathematics \cite{krause2015positive,gau04}. This particular notion of monotonicity should not be confused with that used in convex analysis \cite{Bauschke2017}, which is a different concept in general.
\end{remark}

Note that some references explicitly add positivity to the definition of SI mappings, but we recall that this requirement can be omitted because of the next result.

\begin{fact}\cite{Leung2004}
Let $f\colon\sgens{k}_+\to\sgens{k}_+$ be an SI mapping in the sense of Definition~\ref{def.si}. Then $f(\sgens{k}_+)\subseteq\Intr.$
\end{fact}

Before even trying to compute the fixed point of a SI mapping, we should first answer the more fundamental question of whether the mapping has a fixed point. The result shown in Fact~\ref{sr1} below can be used for this purpose, and note that Fact~\ref{sr1} also establishes uniqueness of the fixed point. To understand the statement, we first need to introduce the concepts of asymptotic mappings and (nonlinear) spectral radius, defined as follows:

\begin{definition}
	\cite{Cavalcante2019_IEEE_TSP} Let $f\colon\sgens{k}_+\to\Intr$ be a continuous SI mapping. The asymptotic mapping associated with $f$ is the continuous mapping defined by
	\begin{equation}
		f_\infty:\sgens{k}_+\to\sgens{k}_+:x\mapsto \lim_{p\to\infty}\dfrac{1}{p}f(p x).
	\end{equation}
	We recall that the above limit always exists in any normed vector space.
\end{definition}


\begin{definition} \label{def.spec_rad}
	The (nonlinear) spectral radius of a continuous SI mapping $f\colon\sgens{k}_+\to\Intr$ is given by the largest eigenvalue of the corresponding asymptotic mapping \cite{Oshime1992,Cavalcante2019,Cavalcante2019_IEEE_TSP}:
	\begin{multline*}
		\rho(f_\infty)=\\\max\{\lambda\in\mathbb{R}_+~|~\ \exists x\in\sgens{k}_+\backslash\{0\}\textrm{ s.t. }f_\infty(x)=\lambda x\}\in\mathbb{R}_+.
	\end{multline*}
	The existence of such an eigenvalue is established in \cite{Oshime1992}, and we remark that there exist simple iterative methods to compute $\rho(f_\infty)$ \cite[Remarks~2 and 3]{Cavalcante2019_IEEE_TSP}.
\end{definition}

\begin{fact} \label{sr1}
	If $f\colon\sgens{k}_+\to\Intr$ is a continuous SI mapping, then $f$ has a fixed point if and only if $\rho(f_\infty)<1$ \cite{Cavalcante2019_IEEE_TSP}. Furthermore, if a fixed point exists, then it is unique \cite{yates95}.
\end{fact}

In many studies, convergence of fixed point iterations is proved using arguments valid for any SI mapping \cite{yates95,huang1998rate,Cavalcante2016,Cavalcante2019_IEEE_TSP,martin11,slawomir09,you2020note,shindoh2020,cavalcante2016low,ho2014data,piotrowski2021fixed}. Therefore, we are unable to establish fine properties of the fixed point iteration, such as geometric convergence in normed vector spaces, because properties of this type do not hold in general \cite[Example~2]{fey2012}. However, the problems being addressed in many studies typically involve only positive concave mappings, which are described below in Definition~\ref{PC}. We remark that these mappings can be seen as a further restriction imposed on SI mappings, as shown in Fact~\ref{PC-SI}. Nevertheless, positive concave mappings remain general enough to address applications not only in the wireless domain, but in many other fields, so we focus on (continuous) positive concave mappings in this study.   

\begin{definition} \label{PC}
Let $f\colon\sgens{k}_+\to\Intr$ be concave with respect to (w.r.t.) the cone order:
\begin{multline}
\forall x,y\in\sgens{k}_+\quad \forall t\in(0,1)\\ f(tx+(1-t)y)\geq tf(x)+(1-t)f(y).
\end{multline}
Then $f$ is called a positive concave (PC) mapping.
\end{definition}


\begin{fact}\cite{Cavalcante2016,Cavalcante2019_IEEE_TSP} \label{PC-SI}
Let $f\colon\sgens{k}_+\to\Intr$ be a PC mapping. Then $f$ is an SI mapping (note: the converse does not hold in general).
\end{fact}

To analyze the convergence of fixed point iterations of monotone mappings on positive or nonnegative cones, researchers in mathematics have extensively used  the metric space $(\Intr, d_T)$ \cite{Lemmens2012,krause2015positive}, where $d_T$ denotes Thompson's metric defined below. Recalling that PC mappings in Definition \ref{PC} are also monotone and defined on nonnegative cones, we can build upon a huge body of the mathematical literature to prove the main results in the next sections.

\begin{definition}
	\label{tmetric}
We define Thompson's metric by
\begin{equation} 
	\label{eq.thompson} 
	\begin{array}{rl}
d_T:\Intr\times \Intr&\to [0,\infty)\\ 
(x,y)&\mapsto\ln(\max\{M(x,y),M(y,x)\}),
\end{array}
\end{equation}
where $M(x,y)\df\inf\{\beta>0~|~x\leq\beta y\}$ and analogously for $M(y,x).$
\end{definition}
Hereafter, we call the metric space $(\Intr, d_T)$ the \emph{Thompson metric space}. The chief results of this paper address properties of positive concave (PC) mappings on compact sets. Therefore, we need to establish equivalence results for compact sets between Thompson metric space and normed vector spaces. To this end, we first state  the following technical results in Remarks \ref{isometric} and~\ref{isometric2}.

\begin{remark} \label{isometric}
The Thompson metric space $(\Intr, d_T)$ is isometric (and, hence, also homeomorphic) to $(\sgens{k}, \|\cdot\|_\infty)$, with an isometry given by the componentwise natural logarithm; namely
\begin{equation} \label{eq.homeomorphism}
L\colon\Intr\to\sgens{k}:x\mapsto (\ln x[1],\ln x[2],\dots,\ln x[k]),
\end{equation}
where $x=(x[1],x[2],\dots,x[k])\in\Intr$ \cite[Proposition 2.2.1]{Lemmens2012}. This fact implies that $(\Intr, d_T)$ inherits all topological and metric properties of $(\sgens{k}, \|\cdot\|_\infty).$ We also note that the inverse isometry from $(\sgens{k}, \|\cdot\|_\infty)$ to $(\Intr, d_T)$ is given as
\begin{equation} \label{eq.inv_homeomorphism}
L^{-1}\colon\sgens{k}\to\Intr:y\mapsto (\exp y[1],\exp y[2],\dots,\exp y[k]),
\end{equation}
where $y=(y[1],y[2],\dots,y[k])\in\sgens{k}.$
\end{remark}

\begin{remark} \label{isometric2}
From Remark \ref{isometric} we obtain in particular that the Thompson metric space $(\Intr, d_T)$ is complete (i.e., any Cauchy sequence is convergent). Furthermore, $U\subset\Intr$ is a closed, bounded, or compact set in $(\Intr, d_T)$ if and only if $L(U)$ is a closed, bounded, or compact set in $(\sgens{k}, \|\cdot\|_\infty)$, and by the norm equivalence in $\sgens{k}$, if and only if $L(U)$ is a closed, bounded, or compact set in $(\sgens{k}, \|\cdot\|)$ for any norm $\|\cdot\|.$ 
\end{remark}

A potential drawback of working in the metric space $(\Intr, d_T)$, which can be observed from Remark \ref{isometric2}, is that the topological properties of a set $U\subseteq\Intr$ in Thompson metric space are the same as of $L(U)\subseteq\sgens{k}$ in a normed vector space, but not necessarily of $U$ considered in this normed vector space. For example, given $u\in\Intr$, the set $U\df\{x\in\Intr~|~x\le u\}$ is \emph{closed and unbounded} in Thompson metric space $(\Intr, d_T)$, just like the set $L(U)=\{x\in\sgens{k}~|~x\leq L(u)\}$ in a normed vector space $(\sgens{k}, \|\cdot\|).$ On the other hand, $U$ is clearly bounded and neither closed nor open in a normed vector space. However, as Corollary~\ref{isometric4} below shows, as far as compact sets are concerned, we do not have to worry about the space in context [Thompson metric space $(\Intr, d_T)$ or a normed vector space $(\sgens{k}, \|\cdot\|)$].


\begin{corollary} \label{isometric4}
A set $U\subset\Intr$ is compact in $(\Intr, d_T)$ if and only if $U$ is compact in $(\sgens{k}, \|\cdot\|)$ for any norm $\|\cdot\|.$
\proof\noindent $(\Rightarrow)$ Let $U\subset\Intr$ be compact in $(\Intr, d_T).$ From Remark \ref{isometric2}, $L(U)$ is compact in $(\sgens{k}, \|\cdot\|)$ for any norm $\|\cdot\|.$ As $L^{-1}$ is continuous, $L^{-1}(L(U))=U$ is also compact in $(\sgens{k}, \|\cdot\|).$

\noindent $(\Leftarrow)$ Let $U\subset\Intr$ be compact in $(\sgens{k}, \|\cdot\|)$ for any norm $\|\cdot\|.$ As $L^{-1}$ is an isometry from $(\sgens{k}, \|\cdot\|_\infty)$ to $(\Intr, d_T)$, an analogous argument to that used in Remark \ref{isometric2} allows us to conclude that $L^{-1}(U)$ is compact in $(\Intr, d_T).$ As $L$ is continuous, $L(L^{-1}(U))=U$ is also compact in $(\Intr, d_T).$ \IEEEQEDhere
\end{corollary}

The following corollary establishes further equivalent properties of compact sets in Thompson metric space and in normed vector spaces. More precisely, it provides an equivalent definition of a compact set in Thompson metric space, and note that this definition corresponds to the well known Heine-Borel theorem for finite-dimensional normed vector spaces.

\begin{corollary} \label{isometric3}
A set $U\subset\Intr$ is compact in $(\Intr, d_T)$ if and only if $U$ is closed and bounded in $(\Intr, d_T).$
\end{corollary}
\proof\noindent $(\Rightarrow)$ Let $U\subset\Intr$ be compact in $(\Intr, d_T).$ In view of Remark \ref{isometric2}, $L(U)$ is compact in any normed vector space $(\sgens{k}, \|\cdot\|).$ By the Heine-Borel Theorem, $L(U)$ is closed and bounded in $(\sgens{k}, \|\cdot\|)$. From Remark \ref{isometric2}, $U$ is closed and bounded in $(\Intr, d_T).$

\noindent $(\Leftarrow)$ Let now $U\subset\Intr$ be a closed and bounded set in $(\Intr, d_T).$ In view of Remark \ref{isometric2}, $L(U)$ is closed and bounded in any normed vector space $(\sgens{k}, \|\cdot\|).$ By the Heine-Borel Theorem, $L(U)$ is compact in $(\sgens{k}, \|\cdot\|).$ From Remark~\ref{isometric2}, $U$ is compact in $(\Intr, d_T).$ \IEEEQEDhere

We finish this section with a list of standard definitions and results in fixed point theory.

\begin{definition}
Let $f\colon\Intr\to\Intr.$ We say that $f$ is a para-contraction (is strictly nonexpansive) w.r.t. Thompson's metric on $\Intr$ if
\begin{equation}
\forall x,y\in\Intr\quad x\neq y\implies d_T(f(x),f(y))<d_T(x,y).
\end{equation}
Moreover, we say that $f$ is a $c$-Lipschitz contraction w.r.t. Thompson's metric on $\Intr$ if
\begin{equation} \label{conv1}
\exists c\in[0,1)\quad \forall x,y\in\Intr\quad d_T(f(x),f(y))\leq c\,d_T(x,y).
\end{equation}
We call $c$ in (\ref{conv1}) a contraction factor of $f.$ Finally, we say that $f$ is a local $c$-Lipschitz contraction w.r.t. Thompson's metric on a compact set $U\subset\Intr$ if
\begin{equation} \label{conv2}
\exists c\in[0,1)\quad \forall x,y\in U\quad d_T(f(x),f(y))\leq c\,d_T(x,y).
\end{equation}
We call $c$ in (\ref{conv2}) a local contraction factor of $f$ on $U.$
\end{definition}

\begin{remark} \label{bfpt}
Let $f\colon\Intr\to\Intr$ be a $c$-Lipschitz contraction w.r.t. Thompson's metric on $\Intr.$ Then, from the Banach fixed point theorem, there exists a unique fixed point $x^\star\in\Intr$ of $f$ such that the fixed point iteration in (\ref{fpseq}) converges both linearly and geometrically to $x^\star$ with a (contraction) factor $c$, and, in the latter case, with $\gamma= d_T(x_2,x_1)(1-c)^{-1}.$
\end{remark}

\begin{fact}\cite[Lemma~2.1.7]{Lemmens2012} \label{f1}
Let $f\colon\sgens{k}_+\to\Intr$ be an SI mapping and let $f\restr{\Intr}$ be the restriction of $f$ to $\Intr.$ Then $f\restr{\Intr}$ is a para-contraction w.r.t. Thompson's metric on $\Intr.$ 
\end{fact}



\section{Convergence of the fixed point iteration of continuous PC mappings} \label{main}
The main objective of this study is to prove geometric convergence of the fixed point iteration of an arbitrary continuous PC mapping with nonempty fixed point set, and we recall that the fixed point set of these mappings is a singleton as a consequence of Facts~\ref{sr1} and \ref{PC-SI}. To establish geometric convergence, we start by proving that PC mappings $f\colon\sgens{k}_+\to\Intr$ with $\mathrm{Fix}(f)\neq\emptyset$ are local contractions w.r.t. Thompson's metric $d_T$ in (\refeq{eq.thompson}) in a neighborhood $U$ of the fixed point $x^\star\in\mathrm{Fix}(f)$. Therefore, since $(U, d_T)$ is a complete metric space, then standard arguments based on the Banach fixed point theorem guarantees geometric convergence of the fixed point iteration in $(U, d_T)$ for $x_1\in U$. These results are the main topic of Sect.~\ref{contractivity}. 

In engineering applications such as those in the wireless domain, we are typically not interested in geometric convergence of the fixed point algorithm in the metric space $(U,~ d_T)$ for an unspecified neighborhood $U$ of the fixed point, but in a normed vector space $(\sgens{k},~\|\cdot\|)$ (the choice of the norm is irrelevant because of the equivalence of norms in finite dimensional spaces). Unfortunately, the homeomorphism  in (\ref{eq.homeomorphism}) between $(\Int{\sgens{k}_+},~ d_T)$ and $(\sgens{k},~\|\cdot\|_\infty)$ is nonlinear, so it may be unclear whether the results in Sect.~\ref{contractivity} imply geometric convergence in a normed vector space. This problem is addressed in Sect.~\ref{sect.geometric}, where we prove that geometric convergence of the fixed point iteration in $(U,~ d_T)$ with a factor $c\in~[0,1)$ implies geometric convergence in $(\sgens{k},~\|\cdot\|)$ with the same factor $c$ for any starting point $x_1\in\real_+^k$. 

Since the contraction factor $c$ gives useful information about the convergence speed of the fixed point iteration [see (\ref{geom})], in Sect.~\ref{sect.specrad} we show that the concept of (nonlinear) spectral radius is useful not only for establishing the existence of a fixed point as shown in Fact~\ref{sr1}, but also for bounding $c$. 

\subsection{Local $c$-Lipschitz contractivity}
\label{contractivity}

Before delving into the convergence of the fixed point iteration of PC mappings, we first derive properties required for a mapping to be a contraction w.r.t. Thompson's metric. We begin with the following lemma, which is an extension of \cite[Lemma 2.1.7]{Lemmens2012} to $c\in[0,1).$

\begin{lemma} \label{l1}
Let $f\colon\Intr\to\Intr$ be monotonic. Then $f$ is a $c$-Lipschitz contraction w.r.t. Thompson's metric if and only if
\begin{multline} \label{ssc_l1}
\exists c\in[0,1)\quad\textrm{s.t.}\quad\forall x\in\Intr\quad \forall \lambda>1\\ f(\lambda x)\leq\lambda^cf(x).
\end{multline}

\end{lemma}
\proof\noindent $(\Leftarrow)$ Let $x,y\in\Intr$, so that $d_T(x,y)=\ln\lambda$ for some $\lambda\geq 1.$ If $x=y$, then $d_T(x,y)=0$ and $\lambda=1$, in which case $d_T(f(x),f(y))\leq c\,d_T(x,y)$ for any $c\in[0,1).$ Therefore, we assume below that $\lambda>1.$

\emph{The case $\lambda>1$ ($x\neq y$)}: From the definition of $d_T$, one has that $x\leq\lambda y$ and $y\leq\lambda x.$ From monotonicity of $f$ we also obtain that  $f(x)\leq f(\lambda y)$ and $f(y)\leq f(\lambda x)$, and from (\ref{ssc_l1}) one has that $f(\lambda y)\leq\lambda^cf(y)$ and $f(\lambda x)\leq\lambda^cf(x)$ for some $c\in[0,1).$ As a result,
\begin{equation*}
f(x)\leq f(\lambda y)\leq \lambda^cf(y)\textrm{ and }f(y)\leq f(\lambda x)\leq \lambda^cf(x).
\end{equation*}
Thus,
\begin{multline*}
\beta_1\df \inf\big\{\beta>0~|~f(x)\leq\beta f(y)\big\}\leq\lambda^c\textrm{ and }\\ \beta_2\df \inf\big\{\beta>0~|~f(y)\leq\beta f(x)\big\}\leq\lambda^c,
\end{multline*}
and hence also $\max\{\beta_1,\beta_2\}\leq\lambda^c.$ We now conclude that
\begin{multline*}
d_T(f(x),f(y))=\ln(\max\{\beta_1,\beta_2\})\leq\ln\lambda^c=\\ c\ln\lambda=c\,d_T(x,y).\ 
\end{multline*}
\noindent $(\Rightarrow)$ Conversely, assume that $f$ is a $c$-Lipschitz contraction w.r.t. Thompson's metric, and let $x\in\Intr$ and $\lambda>1$ be chosen arbitrarily. We note that $d_T(x,\lambda x)=\ln\lambda$, and thus
\begin{align*}
d_T(f(x),f(\lambda x))\le c\,d_T(x,\lambda x)=c\ln\lambda=\ln\lambda^c.
\end{align*} 

Let
\begin{multline*}
\beta_3\df \inf\big\{\beta>0~|~f(x)\leq\beta f(\lambda x)\big\}\textrm{ and }\\ \beta_4\df \inf\big\{\beta>0~|~f(\lambda x)\leq\beta f(x)\big\}.
\end{multline*}
From monotonicity of $f$ we have that $f(x)\leq f(\lambda x)$ because $\lambda>1$, thus $\beta_3\leq\beta_4.$ Therefore,
\begin{align} \label{conc_l1}
d_T(f(x),f(\lambda x))=\ln\beta_4\leq\ln\lambda^c.
\end{align}
From the definition of $\beta_4$, we obtain that $f(\lambda x)\leq\beta_4 f(x)$, and from (\ref{conc_l1}) we obtain in particular that $\beta_4\leq\lambda^c$. As a result, we also have $\beta_4 f(x)\leq\lambda^c f(x)$, and the proof of the desired result $f(\lambda x)\leq\lambda^c f(x)$ is now complete.  $\blacksquare$

\begin{remark}
Mappings satisfying (\ref{ssc_l1}) are known as $c$-concave mappings in the mathematical literature, where they are considered in the context of partially ordered subsets (such as cones) of generic Banach spaces; see, for example, \cite{Potter1977,Chengbo2005,Liang2006}.
\end{remark}

Continuous PC mappings (and, in particular, positive affine mappings restricted to $\sgens{k}_+$) do not satisfy the properties in Lemma~\ref{l1} in general, even if the mappings have a fixed point. Therefore, these mappings are not \emph{global} contractions in the metric space $(\Int{\sgens{k}_+},~ d_T)$, and thus Lemma~\ref{l1} cannot be directly used with the Banach fixed point theorem to deduce geometric convergence of the fixed point iteration with PC mappings. Nevertheless, the following corollary shows that a weaker condition than (\ref{ssc_l1}) is sufficient for $f$ to be a \emph{local} $c$-Lipschitz contraction w.r.t. Thompson's metric on a compact set $U\subset\Intr$, and in Proposition \ref{l3} we show that this property is satisfied by PC mappings.

\begin{corollary} \label{saved!}
Let $U\subset\Intr$ be a compact set containing a nonempty open set and let $d_0\df\max_{x,y\in U} d_T(x,y)=\ln\lambda_0$ with $\lambda_0>1.$ Let $f\colon\Intr\to\Intr$ be monotonic. Then $f$ is a $c$-Lipschitz contraction w.r.t. Thompson's metric on $U$ if
\begin{equation} \label{ssc}
\exists c\in[0,1)\quad\textrm{s.t.}\quad\forall x\in U\quad \forall \lambda\in(1,\lambda_0]\quad f(\lambda x)\leq\lambda^cf(x).
\end{equation}
\end{corollary}
\proof Let $x,y\in U.$ Then, $d_T(x,y)=\ln\lambda$ for some $1\leq\lambda\leq\lambda_0.$ The proof now follows from the \emph{if} direction $(\Leftarrow)$ of the proof of Lemma \ref{l1}, with $x,y\in\Intr$ considered therein restricted to $x,y\in U$, and $\lambda\geq 1$ considered therein restricted to $\lambda\in[1,\lambda_0].$ $\blacksquare$

Thanks to Lemma \ref{l1} and Corollary \ref{saved!}, we are able to prove the following proposition on local $c$-Lipschitz contractivity of continuous PC mappings.

\begin{proposition} \label{l3}
Let $f\colon\sgens{k}_+\to\Intr$ be a continuous and concave (w.r.t. to the cone order) mapping (i.e., a continuous PC mapping), and let $U\subset\Intr$ be a compact set containing a nonempty open set. Then $f$ is a local $c$-Lipschitz contraction w.r.t. Thompson's metric on $U.$
\end{proposition}
\proof For convenience, we divide the proof into the following parts.

\emph{1) Global properties of $f$}: From concavity of $f$, one has in particular that, for $t\in(0,1)$ and $x\in\sgens{k}_+$, 
\begin{equation*}
f(tx)=f(tx+(1-t)0)\geq tf(x)+(1-t)f(0).
\end{equation*}
Furthermore,
\begin{equation*}
f(x)=f(t(x/t))\geq tf(x/t)+(1-t)f(0).
\end{equation*}  
Let $\lambda\df t^{-1}>1.$ Then, the above inequality can be rewritten as
\begin{equation*}
f(x)\geq \lambda^{-1}f(\lambda x)+(1-\lambda^{-1})f(0),
\end{equation*}  
thus
\begin{equation*}
\lambda f(x)\geq f(\lambda x)+(\lambda-1)f(0),
\end{equation*}
and consequently
\begin{equation} \label{c1}
f(\lambda x)\leq\lambda f(x)+(1-\lambda)f(0),\quad x\in\sgens{k}_+.
\end{equation}

\emph{2) Local properties of $f$ on $U$}: By our assumptions, $f$ is continuous, and $U\subset\Intr$ is closed and bounded (see Corollary \ref{isometric3}), hence $f(U)$ is also closed and bounded. In particular, there exists $\mu>0$ such that $f(0)\geq\mu f(x)$ for all $x\in U.$ We note that we can always select $\mu\in(0,1)$ because concavity and positivity of $f$ on $\sgens{k}_+$ implies its monotonicity on $\sgens{k}_+$ \cite{Cavalcante2016} (see also Fact~\ref{PC-SI}). Noticing that $1-\lambda<0$, we obtain from~(\ref{c1}) that
\begin{equation*}
f(\lambda x)\leq\lambda f(x)+(1-\lambda)\mu f(x),\quad x\in U. 
\end{equation*}
Thus,
\begin{multline} \label{c2}
f(\lambda x)\leq\lambda\Bigg(1+\frac{1-\lambda}{\lambda}\mu\Bigg)f(x)=\\\lambda\Bigg(\frac{(1-\mu)\lambda+\mu}{\lambda}\Bigg)f(x)=\lambda\eta f(x),\quad x\in U,
\end{multline}
where $\eta\df\frac{(1-\mu)\lambda+\mu}{\lambda}>0$, because $\lambda>1$ and $\mu\in(0,1).$

\emph{3) Exponential form of $(\ref{c2})$}: Using the \emph{change of base} rule of logarithms, we note that $\log_\lambda(\lambda\eta)=\frac{\ln(\lambda\eta)}{\ln\lambda}$, and hence $\lambda\eta=\lambda^{\frac{\ln(\lambda\eta)}{\ln\lambda}}.$ Thus, we obtain from (\ref{c2}) that
\begin{equation} \label{logs}
f(\lambda x)\leq\lambda^{\frac{\ln(\lambda\eta)}{\ln\lambda}}f(x)=\lambda^{c(\lambda)}f(x),\quad x\in U,
\end{equation}
where
\begin{equation*}
c(\lambda)\df\frac{\ln(\lambda\eta)}{\ln\lambda}=\frac{\ln((1-\mu)\lambda+\mu)}{\ln\lambda}
\end{equation*}
is a continuous function of $\lambda>1.$

\emph{4) Monotonic increase of $c(\lambda)$}: We show now that $c(\lambda)$ is monotonically increasing for $\lambda>1.$ Denote
\begin{equation*}
s\df(1-\mu)\lambda+\mu,\quad \lambda>1,
\end{equation*}
and observe that $1<s<\lambda$ for any $\mu\in(0,1).$ Thus,
\begin{equation*}
c(\lambda)=\frac{\ln s}{\ln\lambda}<1,\quad \lambda>1.
\end{equation*}  
Moreover, we have
\begin{equation} \label{cprime}
c'(\lambda)=\frac{\frac{(1-\mu)\ln\lambda}{s}-\frac{\ln s}{\lambda}}{(\ln\lambda)^2},\quad \lambda>1.
\end{equation}  
Define
\begin{equation*}
c_2(\lambda)\df\frac{(\ln s)'}{(\ln\lambda)'}=\frac{(1-\mu)\lambda}{s},\quad \lambda>1.
\end{equation*}
Then $c'_2(\lambda)=\frac{(1-\mu)\mu}{s^2}>0$, hence $c_2(\lambda)$ is increasing for $\lambda>1.$

We note that both $\lambda\mapsto\ln s$ and $\lambda\mapsto\ln\lambda$ are continuous for $\lambda\geq 1$ and differentiable for $\lambda>1$, with $1\mapsto 0$ in both cases. Therefore, using Cauchy's (extended) mean value theorem, we conclude that for any $\lambda_1>1$ there exists $\lambda_2$ such that $1<\lambda_2<\lambda_1$ with $c_2(\lambda_2)=c(\lambda_1).$ Thus, as $c_2(\lambda)$ is increasing for $\lambda>1$, we obtain $c_2(\lambda_1)>c(\lambda_1)$ for $\lambda_1>1.$ Therefore, one has
\begin{equation*}
c_2(\lambda)=\frac{(1-\mu)\lambda}{s}>\frac{\ln s}{\ln\lambda}=c(\lambda),\quad \lambda>1,
\end{equation*}
thus also $\frac{(1-\mu)\ln\lambda}{s}-\frac{\ln s}{\lambda}>0$ for $\lambda>1.$ In view of (\ref{cprime}), we conclude therefore that $c(\lambda)$ is increasing for $\lambda>1$ and any $\mu\in(0,1).$  We also note that
\begin{equation*}
1-\mu<c(\lambda)<1,
\end{equation*}
as $\lim_{\lambda\to 1^+}c(\lambda)=1-\mu$ and $\lim_{\lambda\to\infty}c(\lambda)=1.$

\emph{5) Explicit local contraction factor}: Let now $\lambda_0>1$ satisfy $\max_{x,y\in U} d_T(x,y)=\ln\lambda_0$, and let
\begin{equation*}
c\df c(\lambda_0)<1. 
\end{equation*}
Then, from (\ref{logs}) and the fact that $c(\lambda)\leq c$ for all $1<\lambda\leq\lambda_0$, one has
\begin{equation}
f(\lambda x)\leq\lambda^{c(\lambda)} f(x)\leq\lambda^c f(x),\quad x\in U,\ 1<\lambda\leq\lambda_0.
\end{equation}
Therefore, we conclude that $f$ satisfies inequality (\ref{ssc}) in Corollary~\ref{saved!}, implying that $f$ is a local $c$-Lipschitz contraction w.r.t. Thompson's metric on $U.$ $\blacksquare$ 

Interesting consequences of  Proposition \ref{l3} are discussed below.

\begin{remark}
If $f$ has a fixed point $x^\star\in\Intr$, then $U$ can be selected as a compact neighborhood of $x^\star.$ In this case, Proposition \ref{l3} provides a local contraction factor for $f$ on $U.$ We also note that the convergence of a fixed point iteration of $f$ on $U$ is linear w.r.t. Thompson's metric in this case. 
\end{remark}

\begin{example} \label{ex11}
Let $k=1$ and define $f_1(x)=(1/2)x+1/2$ for $x\in\mathbb{R}_+.$ Let $U=[1/2,3/2].$ Then, $(U,d_T)$ is a complete metric space with $f_1(U)\subseteq U$, and, in Proposition \ref{l3}, we may set $\mu=1/3$ as $f_1(0)=1/2$ and $f_1(3/2)=5/4.$ Moreover, we have $\lambda_0=3$, thus $c\df c(\lambda_0)\approx 0.771$, which is a local contraction factor of $f_1$ on $U.$ Thus, by the Banach fixed point theorem, $f_1$ admits a unique fixed point over $U$, which is $x^\star=1$ in this case.
\end{example}  

\begin{example} \label{ex12}
Let $k=1$ and define $f_2(x)=x+1$ for $x\in\mathbb{R}_+.$ As above, let $U=[1/2,3/2]$, so that $(U,d_T)$ is a complete metric space. Then, in Proposition~\ref{l3}, we may choose $\mu=1/3$ because $f_2(0)=1$ and $f_2(3/2)=5/2$, and, as above, we have $\lambda_0=3$, thus $c\df c(\lambda_0)\approx 0.771$, which is a local contraction factor of $f_2$ on $U.$ We note that we have the same values of $\mu,\lambda_0$, and $c$ as in the preceding example. However, the Banach fixed point theorem cannot be applied in this case because $f_2(U)\subsetneq U.$ Indeed, it is seen that $f_2$ does not have a fixed point over~$U.$
\end{example}  

\subsection{Geometric convergence of fixed point iteration in normed vector spaces}
\label{sect.geometric}

Having established linear, and, hence, geometric convergence of the fixed point iteration with a continuous PC mapping in the metric space $(U, d_T)$, where $U$ is a neighborhood of the fixed point (which is assumed to exist), we now proceed to prove that we also have geometric convergence in any normed vector space $(\mathbb{R}^k, \|\cdot\|)$ for any starting point in $\mathbb{R}^k_+$, and the convergence factors in $(U, d_T)$ and  $(\mathbb{R}^k, \|\cdot\|)$ can be selected to be the same. To this end, we need the next result.  

\begin{proposition} \label{conv}
Let $f\colon\Intr\to\Intr$ be a $c$-Lipschitz contraction w.r.t. Thompson's metric and let $x^\star\in\Intr$ be the unique fixed point of $f.$ Then, the fixed point iteration of $f$ with $x_1\in\Intr$ converges geometrically to $x^\star$ with a factor $c\in[0,1)$ w.r.t. any metric induced by a norm in $\sgens{k}$.  
\end{proposition}
\proof Let $x_1\in\Intr$ and let $(x_n)_{n\in\mathbb{N}}$ be the sequence defined through (\ref{fpseq}). We first note 
that, from Remark \ref{bfpt}, $(x_n)_{n\in\mathbb{N}}$ is bounded in $(\Intr, d_T).$ Thus, from Corollary \ref{isometric4} we conclude that the sequence $(x_n)_{n\in\mathbb{N}}$ is also bounded in $(\sgens{k}, \|\cdot\|)$ for any norm $\|\cdot\|.$ Let $b>0$ be such that $\forall n\in\mathbb{N}\ \|x_n\|\le b$ for a given norm $\|\cdot\|.$ We also have $\|x^\star\|\le b$ because $(x_n)_{n\in\mathbb{N}}$ converges to $x^\star.$ Now, let $\delta>0$ be the normality constant of $\sgens{k}_+$ for this norm (see Definition \ref{nc} in the Appendix). Then, using the property that $\forall n\in\mathbb{N}\ d_T(x_{n+1},x^\star)\le c^n\,d_T(x_1,x^\star)$ (as $f$ is a $c$-Lipschitz contraction w.r.t. Thompson's metric) and Fact \ref{L251} in the Appendix \ref{kru}, we deduce

\begin{multline} \label{co1}
\forall n\in\mathbb{N}\quad \|x_{n+1} - x^\star \| \le b(1+2\delta)(e^{d_T(x_{n+1},x^\star)}-1) \\
\le b(1+2\delta)(e^{c^n d_T(x_1,x^\star)}-1)\\= b(1+2\delta){(e^{c^n d_T(x_1,x^\star)})}
\left(1-\dfrac{1}{e^{c^n d_T(x_1,x^\star)}}\right).
\end{multline}
Then, noticing that $\forall n\in\mathbb{N}$ one has $c^n\in[0,1)$, we obtain $e^{c^n d_T(x_1,x^\star)}\leq e^{d_T(x_1,x^\star)}$, and, consequently,
\begin{equation} \label{co1.5}
\forall n\in\mathbb{N}\quad \|x_{n+1} - x^\star \| \le s\left(1-\dfrac{1}{e^{c^n d_T(x_1,x^\star)}}\right),
\end{equation}
where $s\df b(1+2\delta){e^{d_T(x_1,x^\star)}}.$ Using the well-known bound $1-1/v\le \ln v$ for $v>0$, we obtain
\begin{multline} \label{co2}
\forall n\in\mathbb{N}\quad \|x_{n+1} - x^\star \| \le s \ln {e^{c^n d_T(x_1,x^\star)}} =\\ sc^nd_T(x_1,x^\star)=\gamma c^n,
\end{multline}
where $\gamma\df s\,d_T(x_1,x^\star)$, and the proof is complete. \IEEEQEDhere

We now have all ingredients to establish geometric convergence of the fixed point iteration of PC mappings:
\begin{proposition} \label{p1}
Let $f\colon\sgens{k}_+\to\Intr$ be a continuous and concave (w.r.t. to the cone order) mapping (i.e., a continuous PC mapping) with a fixed point $x^\star\in\Intr.$ Then, the fixed point iteration of $f$ with $x_1\in\real_+^k$ converges geometrically to $x^\star$ with a factor $c\in[0,1)$ w.r.t. any metric induced by a norm in $\sgens{k}.$
\end{proposition}

\proof Without any loss of generality, we can assume that $x_1\in\Intr$ because, if $x_1\in\real_+^k\backslash\Intr$, then $f(x_1)\in\Intr$, so we can consider $f(x_1)$ as the starting point of the iterations. We divide the proof into two parts. We first prove geometric convergence of the fixed point iteration of $f$ in a compact neighborhood $U$ of $x^\star$, and then we extend this property to an arbitrary starting point $x_1\in\Intr.$

\emph{1) Geometric convergence on $U$}: We first note that the fixed point $x^\star$ of $f$ satisfies $x^\star\geq f(x^\star)$, thus $f$ is feasible in the sense of \cite{yates95}. Therefore, by Fact \ref{yates} in the Appendix \ref{kru}, the fixed point iteration of $f$ with $x_1\in\Intr$ produces a sequence $(x_n)_{n\in\mathbb{N}}$ converging to $x^\star$ w.r.t. any norm on $\sgens{k}_+.$ Let $U\subset\Intr$ be a compact neighborhood of $x^\star.$ Then, we can find $n_0\in\mathbb{N}$ such that $\forall n>n_0$ one has $x_n\in U.$ Thus, let $(y_n)_{n\in\mathbb{N}}\subset U$ be a sequence defined as $y_n\df x_{n+n_0}$ for $n\in\mathbb{N}.$ Let us also note that $f$ is a local $c$-Lipschitz contraction on $U$ w.r.t. Thompson's metric in view of Proposition \ref{l3}. Moreover, using the isometry provided in (\ref{eq.homeomorphism}) in Remark \ref{isometric}, we conclude that $(U,d_T)$ is a complete metric space isometric to a complete metric subspace of $(\sgens{k},\|\cdot\|_\infty)$~\cite{BabyRudin}. Therefore, as $(y_n)_{n\in\mathbb{N}}\subset U$ is a fixed point iteration of $f$ converging to $x^\star$ in a complete metric space $(U,d_T)$, and $f$ is a local $c$-Lipschitz contraction on $U$ w.r.t. Thompson's metric, we conclude that $(y_n)_{n\in\mathbb{N}}$ is obviously bounded and that it converges linearly (and hence, also geometrically) to $x^\star$ w.r.t. Thompson's metric. Owing to this fact, we may reproduce now the proof of Proposition~\ref{conv}, with the sequence $(x_n)_{n\in\mathbb{N}}\subset\Intr$ considered therein replaced with the sequence $(y_n)_{n\in\mathbb{N}}\subset U$, and the complete metric space $(\Intr,d_T)$ considered therein replaced with the complete metric space $(U,d_T)$, to conclude that $(y_n)_{n\in\mathbb{N}}$ converges geometrically to $x^\star$ with a factor $c\in[0,1)$ w.r.t. any metric induced by a norm in $\sgens{k}$, with a constant $\gamma_y>0$ determined through (\ref{co1})-(\ref{co2}).

\emph{2) Geometric convergence on $\Intr$}: If $n_0>1$, then we need to show that the sequence $(x_n)_{n\in\mathbb{N}}$ also has this property. To this end, let us define $\eta\df\max_{n\in\{1,\dots,n_0\}}\|x_{n+1}-y_{n+1}\|$, and assume that $c=0.$ Then, $f$ is constant on $U$, hence $\forall n\in\mathbb{N}$ one has $y_n=x_{n_0+n}=x^\star.$ Thus, taking into consideration that for any $v\in(0,1)$, we obtain

\begin{equation*}
\forall n\in\{1,\dots,n_0\}\quad \|x_{n+1}-x^\star\|\leq\eta=\Big(\frac{\eta}{v^n}\Big)v^n\leq\Big(\frac{\eta}{v^{n_0}}\Big)v^n,
\end{equation*}
and we conclude that the sequence $(x_n)_{n\in\mathbb{N}}$ converges geometrically to $x^\star$ with a factor $v\in(0,1)$ w.r.t. any metric induced by a norm in $\sgens{k}$, with a constant $\eta/v^{n_0}.$ Assume now that $c\in(0,1).$ Then,
\begin{multline*}
\forall n\in\{1,\dots,n_0\}\quad \|x_{n+1}-x^\star\|\leq\\\|x_{n+1}-y_{n+1}\|+\|y_{n+1}-x^\star\|
\leq \eta+\gamma_yc^n\\=\Big(\frac{\eta}{c^n}+\gamma_y\Big)c^n\leq\Big(\frac{\eta}{c^{n_0}}+\gamma_y\Big)c^n=\gamma_{x_1}c^n,
\end{multline*}
where $\gamma_{x_1}\df\eta c^{-n_0}+\gamma_y.$ Moreover, 
\begin{multline*}
\forall n>n_0\quad \|x_{n+1}-x^\star\|=\|y_{n+1-n_0}-x^\star\|\leq\gamma_yc^{n+1-n_0}\\
=c^{1-n_0}\gamma_yc^n=\gamma_{x_2}c^n,
\end{multline*}
where $\gamma_{x_2}\df c^{1-n_0}\gamma_y.$ Let $\gamma_x\df\max\{\gamma_{x_1},\gamma_{x_2}\}.$ We have found $c\in(0,1)$ and $\gamma_x>0$ such that $\forall n\in\mathbb{N}$ one has $\|x_{n+1}-x^\star\|\leq\gamma_x c^n.$ \IEEEQEDhere 

An interesting question is whether the fixed point iteration is guaranteed to converge geometrically (in a normed vector space) with slightly more general classes of concave mappings than those considered in Proposition~\ref{p1}. The next example shows that this special property of PC mappings  does not necessarily hold even if we simply replace the assumption of positivity by nonnegativity. 

\begin{figure}
	\begin{subfigure}{\linewidth}
		\begin{center}
			\includegraphics[width=.9\columnwidth]{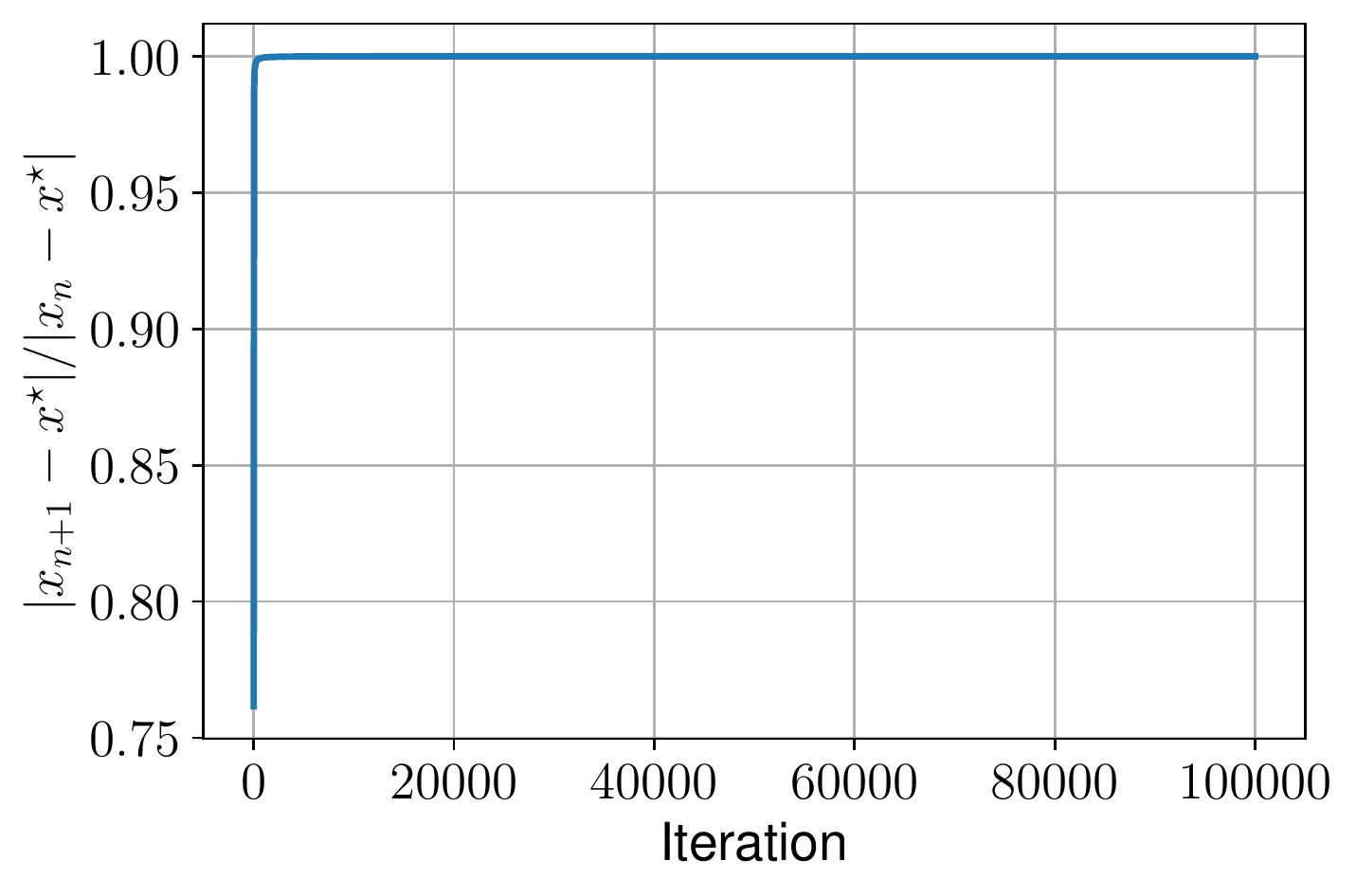}
			\caption{}
		\end{center}
    \end{subfigure}

	\begin{subfigure}{\linewidth}
		\begin{center}
			\includegraphics[width=.9\columnwidth]{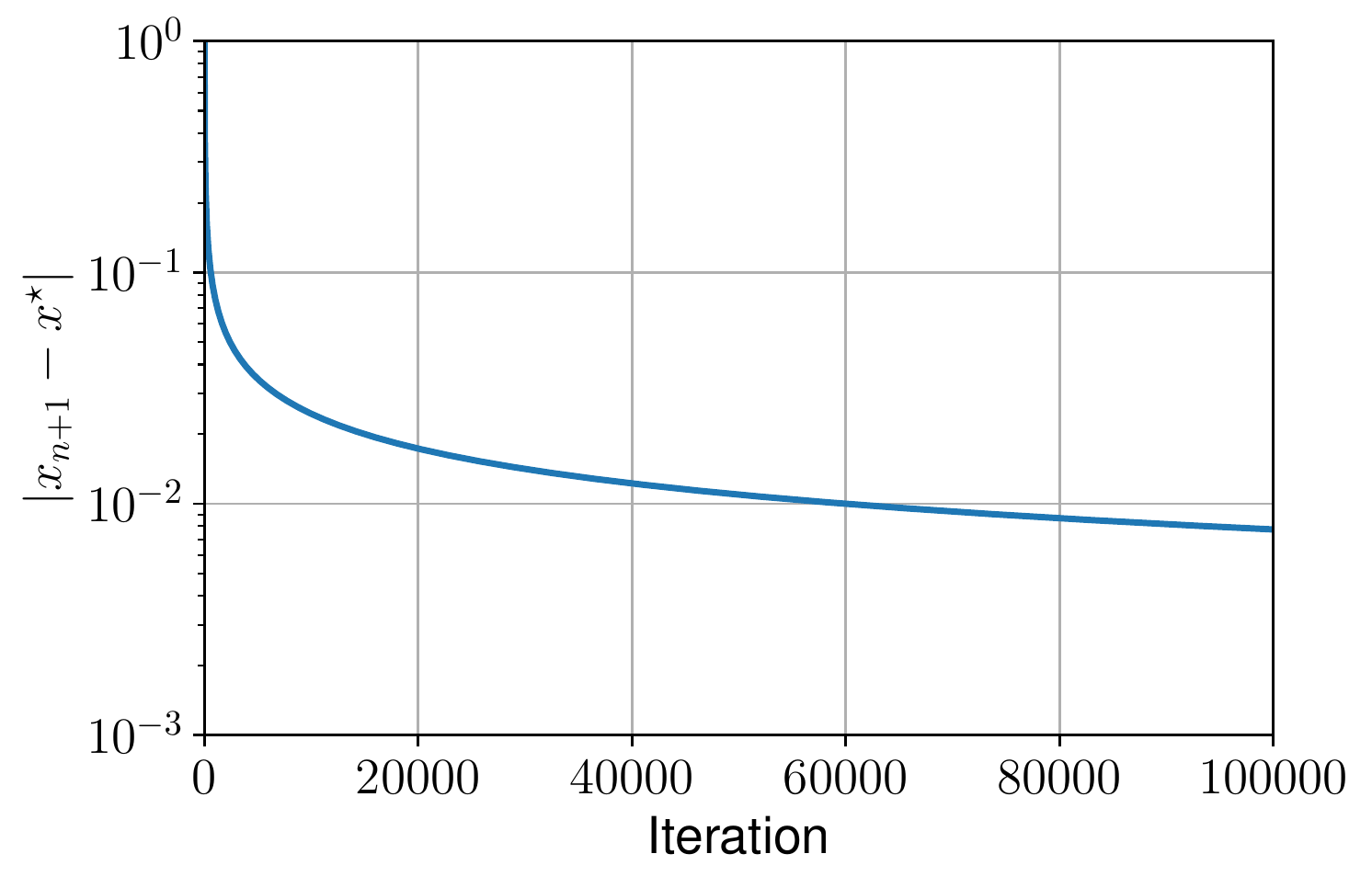}
			\caption{}
		\end{center}
	\end{subfigure}
	\caption{Fixed point iteration of the mapping in Example~\ref{renato_ex}. In this figure $x^\star\in\mathrm{Fix}(g)$ is the limit of the sequence $(x_n)_{n\in\mathbb{N}}$, and note that the fixed point is not unique. (a) Convergence rate. (b) Estimation error.}
    	\label{fig.sublinear}		
\end{figure}

\begin{figure}
	\begin{subfigure}{\linewidth}
		\begin{center}
			\includegraphics[width=.9\columnwidth]{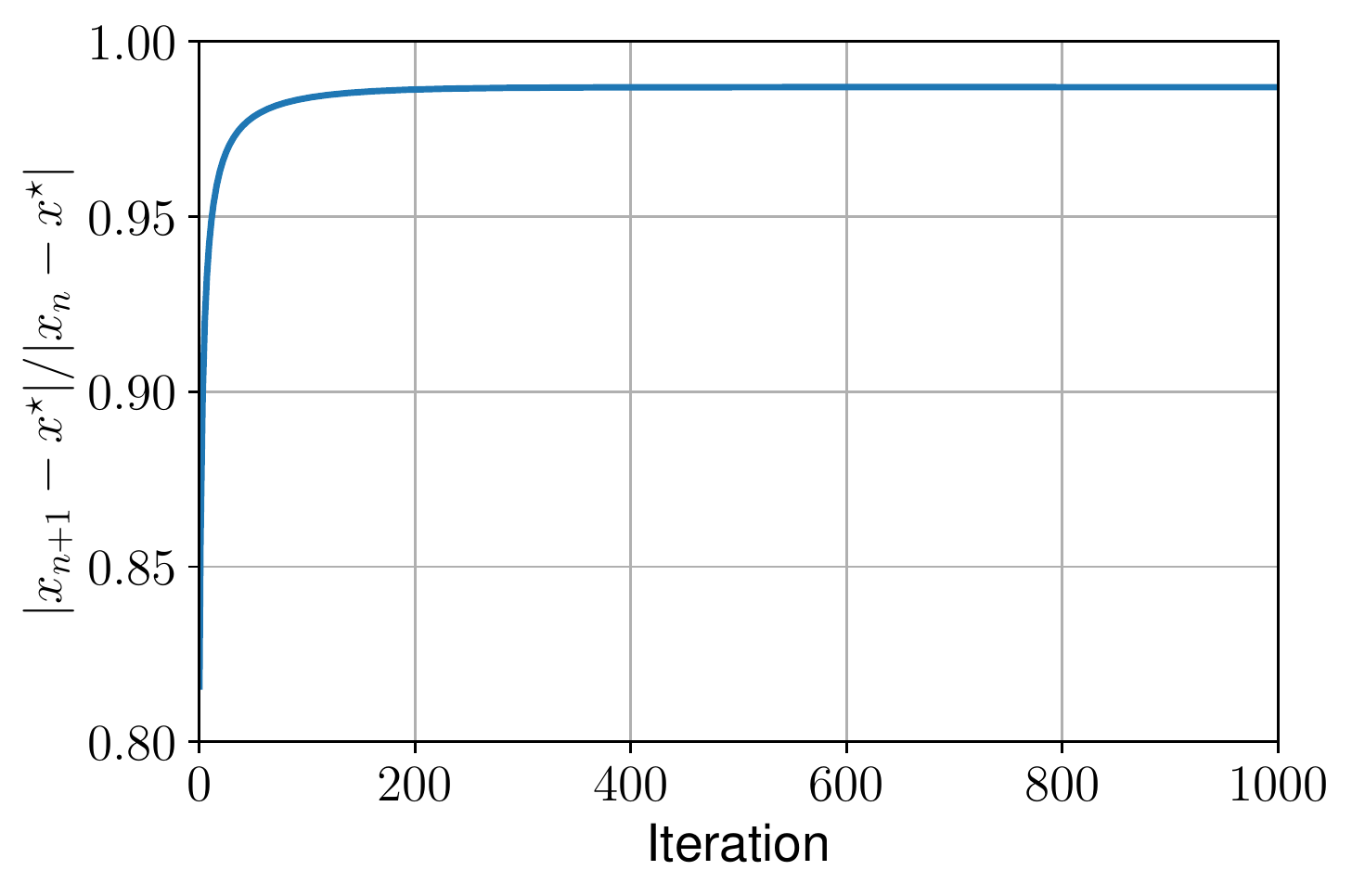}
			\caption{}
		\end{center}
	\end{subfigure}
	
	\begin{subfigure}{\linewidth}
		\begin{center}
			\includegraphics[width=.9\columnwidth]{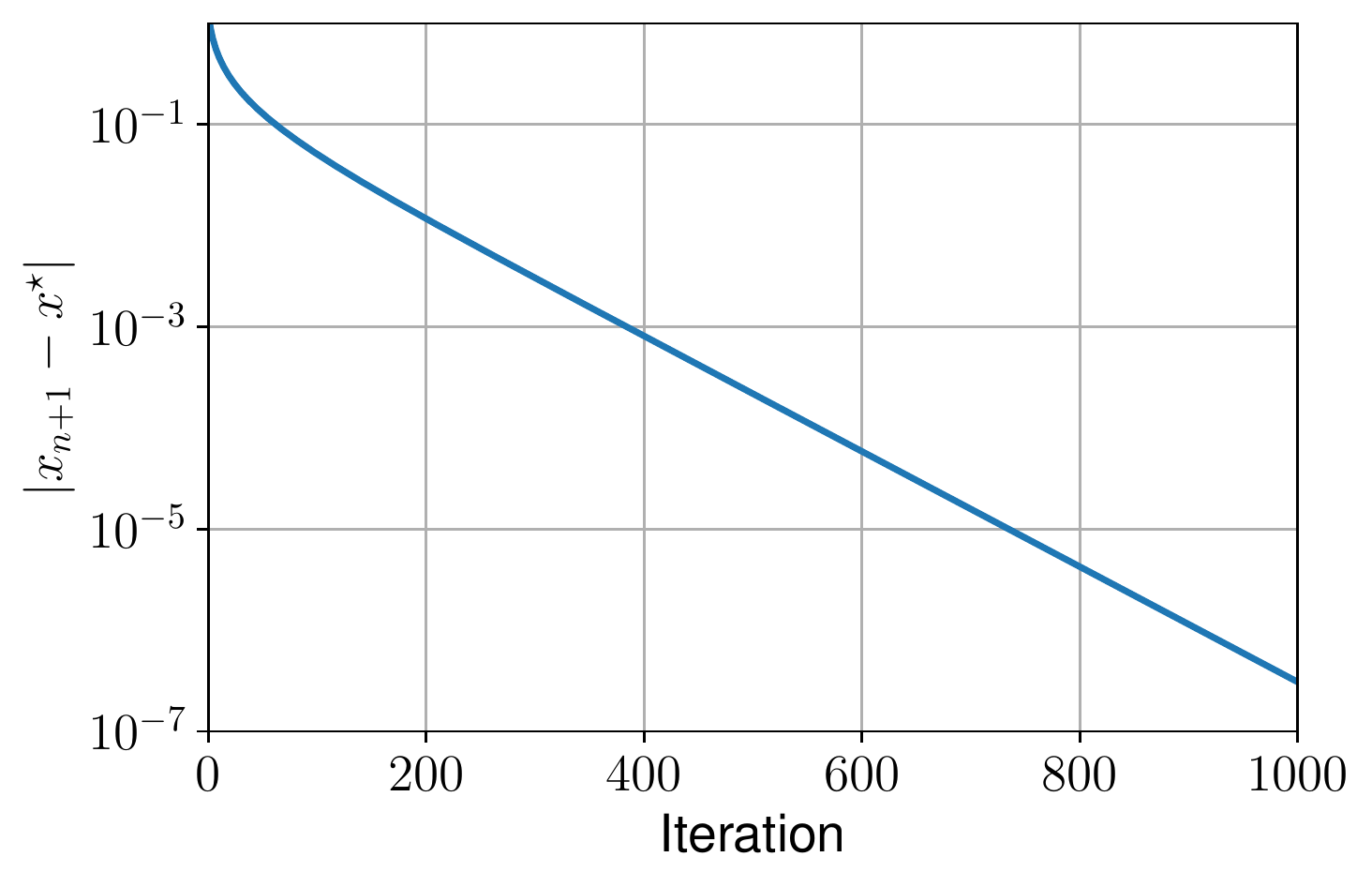}
			\caption{}
		\end{center}
	\end{subfigure}
		\caption{Fixed point iteration of the mapping in Example~\ref{example.linear}. (a) Convergence rate. (b) Estimation error.}
	\label{fig.linear}	
\end{figure}

\begin{example} \label{renato_ex}
In \cite{fey2012}, the authors consider the mapping $$f\colon\mathbb{R}_+\to\Int{\mathbb{R}}:x\mapsto 4/(1+e^{2-x}),$$ and they show that $f$ is an SI mapping with the property that the fixed point iteration of $f$ converges sublinearly to its unique fixed point $x^\star=2.$ In particular, the fixed point iteration of $f$  satisfies neither (\ref{geom}) nor (\ref{linear}) for $\tilde{x}=x^\star$ and $x_{n+1}=f(x_n)$ with $x_1\in\mathbb{R}_+$, $n\in\mathbb{N}.$ This example motivates us to consider the function $g\colon\mathbb{R}_+\to\mathbb{R}_+$ such that $g(x)=x$ for $x\in[0,2]$ and $g(x)=f(x)$ for $x>2.$ Namely, let
\begin{equation}
	\label{eq.g}
  g(x)=\left\{
  \begin{array}{ll}
    x & \textrm{ for }0\leq x\leq 2\\
    \frac{4}{(1+e^{2-x})} & \textrm{ for }x>2.\\
  \end{array}\right.
\end{equation}  
Then, through algebraic manipulations, the following expressions of the first and second derivatives of $g$ can be obtained:
\begin{equation*}
  g'(x)=\left\{
  \begin{array}{ll}
    1 & \textrm{ for }0\leq x\leq 2\\
    \frac{4e^{2+x}}{(e^2+e^x)^2} & \textrm{ for }x>2,\\
  \end{array}\right.
\end{equation*}  
and
\begin{equation*}
  g''(x)=\left\{
  \begin{array}{ll}
    0 & \textrm{ for }0\leq x\leq 2\\
    \frac{4e^{2+x}(e^2-e^x)}{(e^2+e^x)^3} & \textrm{ for }x>2.\\
  \end{array}\right.
\end{equation*}  
We verify that both $g'$ and $g''$ are continuous on $\mathbb{R}_+$ and that $g''(x)\leq 0$ for $x\in\mathbb{R}_+$, implying that $g$ is a continuous and concave function. Moreover, if $x_1>2$, then the fixed point iterations of $g$ and $f$ generate the same sequence; i.e., $\forall n\in\mathbb{N}$  $x_{n+1}=f(x_n)=g(x_n)$. Therefore, from \cite[Example~2]{fey2012} we obtain that the fixed point iteration of $g$ converges sublinearly to $x^\star=2$ for $x_1>2$ irrespective of the selection of a compact neighborhood of $x^\star=2$, for example, on the interval $U=[1,x_1].$ Fig.~\ref{fig.sublinear} illustrates the convergence of the fixed point iteration of $g$ starting at $x_1=4$, and we note that the ratio ${|x_{n+1}-x^\star|}/{|x_n-x^\star|}$ converges to the value one, which is the expected behavior of sequences that converge sublinearly. To highlight further the slow convergence of the algorithm, we also remark that the error $|x_{n+1}-x^\star|$ remains greater than $10^{-3}$ with $n=10^5$ iterations in this experiment. Thus, at first sight, this example invalidates in particular the assertion of Proposition~\ref{p1} on geometric convergence of the fixed point iteration of $g.$ However, we note that $g(0)=0$, so none of the assertions in Propositions \ref{l3}-\ref{p1} apply to $g$ because $g$ is not a PC mapping in the sense of Definition~\ref{PC}. Therefore, this example underlines the importance of positivity of the mappings under consideration for the results in the above propositions to hold.  
\end{example}

\begin{example}
	\label{example.linear}
	By adding a small constant $\epsilon>0$ to the function $g$ in (\ref{eq.g}), we obtain a continuous PC mapping with a unique fixed point denoted by $x^\star$, and thus we re-establish geometric convergence of the fixed point iteration as an implication of Proposition~\ref{p1}. This fact is illustrated in Fig.~\ref{fig.linear}, where we set $\epsilon=10^{-3}$. The result in this figure indicates linear convergence of the fixed point iteration, and we recall that linear convergence implies geometric convergence. In this experiment, the error $|x_{n+1}-x^\star|$ is smaller than $10^{-3}$ already with $n=400$ iterations, and with around $2,000$ iterations the error is reaching the machine precision of the computer used for the simulation. 
\end{example}
	
\begin{remark}
	A potential caveat of adding a small positive constant to a nonnegative concave function is that we may lose its fixed points. This fact can be verified for the function $g$ in (\ref{eq.g}) in Example \ref{renato_ex}. Its fixed point set $\mathrm{Fix}(g)$ is the interval $[0, 2]$, but $\mathrm{Fix}(g+\epsilon)$ is only a singleton for any $\epsilon>0.$ 
\end{remark}



\subsection{Spectral radius and the convergence rate}
\label{sect.specrad}
It follows from (\ref{geom}) that the convergence factor $c$ gives an indication on the number of iterations required to obtain a good approximation of the fixed point of a PC mapping. Therefore, the results derived so far should be supplemented with a computationally simple means of obtaining information about~$c$. A result of this type is provided in this section. In particular, we note that Proposition \ref{prop.linearc} proves that, under mild assumptions, the spectral radius of a PC mapping is a lower bound on $c$, and we recall that simple algorithms for computing the spectral radius are available in the literature~\cite[Remarks 2 and 3]{Cavalcante2019_IEEE_TSP}. 

\begin{proposition} \label{prop.linearc}
Let $f\colon\sgens{k}_+\to\Intr$ be a continuous and concave (w.r.t. to the cone order) mapping (i.e., a continuous PC mapping) with a fixed point $x^\star\in\Intr.$ Let $U\subset\Intr$ be a compact neighborhood of $x^\star.$ Furthermore, let $c\in [0,1)$ be a local contraction factor of $f$ on $U$, and denote by $\rho\df\rho(f_\infty)$ the spectral radius of the asymptotic mapping $f_\infty$ of $f$.  Then $c\ge \rho.$
\end{proposition}
\proof Denote by $v\in\sgens{k}_+\backslash\{0\}$ a (nonlinear) eigenvector associated with the spectral radius $\rho$; i.e., a vector $v\in\sgens{k}_+\backslash\{0\}$ satisfying $f_\infty(v)=\rho v$, which is guaranteed to exist \cite{Nussbaum1986}. The result is trivial if $\rho=0$, so we focus on the case $\rho>0$. Let $x^\star\in U$ be the unique fixed point of $f$ (see Fact~\ref{sr1}) and let $(x_n)_{n\in\mathbb{N}}$ be constructed via the fixed point iteration of $f$ with $x_1\in U$ such that $x_1\ll x^\star$ or $x_1\gg x^\star.$ Then, from Proposition \ref{p1}, we can find a constant $\gamma>0$ such that
\begin{align} \label{eq.linearconv}
\forall n\in\mathbb{N}\quad \|x_{n+1}-x^\star\|\leq\gamma c^n.
\end{align}
Furthermore, there exists $\epsilon>0$ such that $x_1\leq x^\star-\epsilon v$ or $x_1\geq x^\star+\epsilon v$, so \cite[Proposition~5(iii)]{cavalcante2018spectral} asserts that
\begin{align}
	\label{eq.lowerbound}
\forall n\in\mathbb{N}\quad \rho^n\epsilon\|v\|\leq\|x_{n+1}-x^\star\|.
\end{align}
By combining the inequalities in (\ref{eq.linearconv}) and (\ref{eq.lowerbound}), we conclude that $\forall n\in\mathbb{N}$ one has that $\rho^n\epsilon\|v\|\leq\gamma c^n$. Therefore $0<\epsilon\|v\|\leq\gamma (c/\rho)^n$ for all $n\in\mathbb{N}.$ If $c<\rho$, then $c/\rho<1$ and $\lim_{n\to\infty}(c/\rho)^n=0.$ As a result, there exists $n_0\in\mathbb{N}$ such that $\epsilon\|v\|>\gamma(c/\rho)^{n_0}.$ This contradiction implies that we must have $c\geq\rho$, and the proof is complete. \IEEEQEDhere


\begin{remark}
Let $f\colon\sgens{k}_+\to\Intr$ be a continuous PC mapping with a fixed point $x^\star\in\Intr.$ Let $U\subset\Intr$ be a compact neighborhood of $x^\star.$ Then, it is important to highlight that the meaning of $c\in[0,1)$ is consistent among
\begin{itemize}
\item Proposition \ref{l3}, where $c$ is used to denote the local contraction factor of $f$ on $U$,
\item Proposition \ref{prop.linearc}, where $c$ is used in the sense of Proposition~\ref{l3}, 
\item Proposition \ref{p1}, where $c$ is used in the sense of Propositions~\ref{l3} and \ref{prop.linearc}, and it is also recognized as a factor of geometric convergence of the fixed point iteration of $f$ to the unique fixed point.
\end{itemize}    
In particular, Part 5) of the proof of Proposition \ref{l3} provides a method to determine the value of $c$ as 
$c\df c(\lambda_0)=\frac{\ln((1-\mu)\lambda_0+\mu)}{\ln\lambda_0}$, with $\lambda_0=e^{\max_{x,y\in U} d_T(x,y)}\in(1,\infty)$ and $\mu\in(0,1)$ such that $f(0)\geq\mu f(x)$ for all $x\in U.$ Therefore, it is seen that, if $f$ is increasing slowly in the sense that $x^\star=f(x^\star)$ is not much larger than $f(0)$, and the neighborhood $U$ of $x^\star$ is sufficiently small, then the value of $\mu$ can be fixed to a value close to 1. Similarly, the smaller the neighbourhood $U$ of $x^\star$ is considered, the closer the value of $\lambda_0$ is to 1, with $\lim_{\lambda\to 1^+}c(\lambda)=1-\mu.$ Thus, we obtain a useful insight into the relationship between $f$ and $U$, and we also obtain an achievable lower bound for $c$ alternative to the one provided in Proposition \ref{prop.linearc}, where $c$ is bounded below by the spectral radius of of the asymptotic mapping $f_\infty$ of $f.$ 
\end{remark}

\section{Applications}
\label{applications}

We now discuss two common applications of the fixed point iteration of positive concave mappings in wireless networks. We also show how the theory developed in the previous sections is able to explain, with rigorous mathematical statements, convergence properties of existing algorithms.

\subsection{Joint optimization of transmit power, user-base station assignment, and receive beamformers in massive MIMO systems}
\label{sect.power_control}


Consider the uplink of a wireless network with $k$ single antenna users sending data to $m$ base stations, each equipped with $L$ antennas. The set of users and base stations are given by $\mathcal{U}=\{1,\ldots,k\}$ and $\mathcal{B}=\{1,\dots,m\}$, respectively. The set of base stations to which user $u\in\mathcal{U}$ is allowed to connect is represented by the set $\mathcal{B}_u\subseteq\mathcal{B}$.  Denote by $x=[x[1],\dots,x[k]]\in\sgens{k}_{+}$ the transmit power vector, where the $u$th coordinate $x[u]\ge 0$ of the vector $x$ is the transmit power of user $u\in\mathcal{U}$. Let $h_{u,b}$, with  $(u,b)\in\mathcal{U}\times\mathcal{B}$, represent the channel between user $u$ and base station $b$. As common in the literature, for every $(u,b)\in\mathcal{U}\times\mathcal{B}$, we assume $h_{u,b}$ to be a random vector taking values in $\mathbb{C}^L$.  For simplicity, we further assume $h_{u,b}$ to be zero mean with positive definite channel covariance matrix $R_{u,b}\df E[h_{u,b}h_{u,b}^H]\in\mathbb{C}^{L\times L}$.  In this setting, the signal-to-interference-noise ratio (SINR) of user $u$ as a function of its receive beamformer $v\in\mathcal{V}_u$, base station assignment $b$, and transmit power vector $x\in\sgens{k}_{+}$ is given by \cite[Sect.~1.4.2]{martin11}
\begin{align*}
	\begin{array}{rl}
		s_{u}:\sgens{k}_+\times\mathcal{V}_u\times\mathcal{B}_u&\to\real_+ \\ 
		(x,v,b)&\mapsto \dfrac{x[u]~v^H R_{u,b}v}{\sum_{j\in\mathcal{U}\backslash\{u\}}v^H(x[j]R_{j,b}+\sigma^2~I_{L})v},
	\end{array}
\end{align*}
where $\sigma^2>0$ denotes the noise power; $I_{L}\in\sgen{L}{L}$ is the identity matrix;  and $\mathcal{V}_u\subset\mathbb{C}^L$ is a set of allowed beamformers for user $u$, typically chosen to be $\mathcal{V}_u=\{v \in \mathbb{C}^L~|~{v}^H {v}=1\}$. As in \cite{hanly1995algorithm}, users are free to choose the base station to be served.

A standard power control problem is to find the transmit power, the user-base station assignments, and the beamformers that minimize the sum power while satisfying given SINR requirements. Formally, denoting by $\gamma_u>0$ the minimum SINR requirement of user $u\in\mathcal{U}$, we pose the optimization problem as follows:

\begin{align}
	\begin{array}{rl}
		\text{minimize}_{x\in\sgens{k}_+}&\sum_{u\in\mathcal{U}} x[u] \\
		\text{subject to} & \forall u\in\mathcal{U}\quad \max_{(v,b)\in\mathcal{V}_u\times\mathcal{B}_u}s_u(x,v,b)\ge \gamma_u.
	\end{array}
	\label{problem.opt}
\end{align}
For later reference, we say that a vector $x\in\real_+^k$ is feasible to problem (\ref{problem.opt}) if it satisfies all constraints, even if $x$ is not necessarily a solution. In the above formulation, the optimization of the beamformers and user-base station assignments is implicit in the constraint. More precisely, if a power allocation $x$ is given (e.g., an optimal power allocation $x^\star$), then a beamformer $v^\star_u$ and a base station assignment $b_u^\star$  maximizing the SINR for user $u\in\mathcal{U}$ can be recovered by computing:
\begin{align*}
	(v^\star_u,b_u^\star)\in\mathrm{argmax}_{(v,b)\in\mathcal{V}_u\times\mathcal{B}_u} s_u(x,v,b).
\end{align*}
In practice, given $x\in\sgens{k}_+$ and user $u\in \mathcal{U}$, the tuple $(v^\star_u,b_u^\star)$ is obtained by, first, computing the beamformers maximizing the SINR for each of the $|\mathcal{B}_u|\le m$ possible assignments, and by, second, checking the choice that leads to the maximum SINR. The computational complexity of this approach is manageable  if for each user $u$ the number $|\mathcal{B}_u|$ of candidate assignments is not too large, which is a typical property of cellular networks because base stations that are too far away are unable to provide adequate SINR, so they can be ruled out from consideration. We also note that, if fully digital beamformers are used, in which case we can choose $\mathcal{V}_u=\{v \in \mathbb{C}^L~|~{v}^H {v}=1\}$ for every $u\in\mathcal{U}$, then the beamformer maximizing the SINR for a given user, base station assignment, and power vector can be easily computed by maximizing a generalized Rayleigh quotient. Nevertheless, we do not necessarily assume digital beamformers in this section; we only require the following very mild assumption:

\begin{assumption}
	\label{assumption.beamformers}
The set $C_{u,x}=\mathrm{argmax}_{(v,b)\in\mathcal{V}_u\times\mathcal{B}_u} s_u(x,v,b)$ is nonempty for every $u\in\mathcal{U}$ and every $x\in\sgens{k}_+$. 
\end{assumption}

A well-known result in the literature shows that, if a vector $x^\star$ is a solution to the optimization problem (\ref{problem.opt}) (assuming that a solution exists), then the constraints have to be satisfied with equality; i.e., 
\begin{align}
	\label{eq.constraint}
	\forall u\in\mathcal{U}~\max_{(v,b)\in\mathcal{V}_u\times\mathcal{B}_u}s_u(x^\star,v,b)= \gamma_u.
\end{align}
This result can be easily proved by noticing that, if a strict SINR inequality is obtained for some user, then we can decrease the transmit power of this user (thus decreasing the cost function in $(\ref{problem.opt})$),  until an equality is achieved. By decreasing the transmit power of the user, the SINR of other users cannot decrease, so the resulting power vector is also feasible. Therefore, feasible power allocations satisfying all constraints with at least one strict inequality cannot be optimal if problem $(\ref{problem.opt})$ has a solution.

The above fact enables us to replace the inequalities in the constraints of problem (\ref{problem.opt}) with equalities. Therefore, using simple algebraic manipulations, we can equivalently express the constraints in (\ref{problem.opt}) as the fixed point set of the following mapping $f:\sgens{k}_+\to\Intr$:
\begin{align}
	\label{equivalent.constraint}
	\forall u\in\mathcal{U}~\max_{(v,b)\in\mathcal{V}_u\times\mathcal{B}_u}s_u(x^\star,v,b)= \gamma_u \iff x^\star\in\mathrm{Fix}(f),
\end{align}
where 
\begin{align}
	\label{eq.int_map_pc}
f:\sgens{k}_+\to\Intr:x\mapsto [f_1(x),\ldots,f_k(x)]
\end{align}
 and $\forall u\in\mathcal{U}$
\begin{align*}
	\begin{array}{rl}
		f_{u}:\sgens{k}_+&\to\Intrs \\ 
		x\mapsto &\min\limits_{(v,b)\in\mathcal{V}_u\times\mathcal{B}_u} \gamma_u\dfrac{\sum_{j\in\mathcal{U}\backslash\{u\}}v^H(x[j]R_{j,b}+\sigma^2~I_{L})v}{v^H R_{u,b}v}.
	\end{array}
\end{align*}
To avoid digressions of little practical relevance, assume that $k\geq 2$ and $\forall u\in\mathcal{U}\quad 0\notin\mathcal{V}_u$, so that $\forall v\in\mathcal{V}_u\quad v^H R_{u,b} v>0$ because every matrix $R_{u,b}$ is positive definite. Since $\gamma_u>0$ and $\sigma^2>0$ by assumption, we deduce $\forall x\in\real_+^k\quad \forall u\in\mathcal{U}$:
\begin{multline}
	\label{eq.bounded_away}
f_u(x) = \min\limits_{(v,b)\in\mathcal{V}_u\times\mathcal{B}_u} \gamma_u\dfrac{\sum_{j\in\mathcal{U}\backslash\{u\}}v^H(x[j]R_{j,b}+\sigma^2~I_{L})v}{v^H R_{u,b}v} \\
\ge \min\limits_{(v,b)\in\mathcal{V}_u\times\mathcal{B}_u}\gamma_u\sigma^2(k-1)\dfrac{v^H v}{v^H R_{u,b}v} = f_u(0) > 0.
\end{multline}
Furthermore, each function $f_u$ can be seen as the minimum of affine functions, so $f_u$ is not only positive but also concave, and thus $f_u$ is a PC mapping. We can tacitly assume $f_u$ to be continuous in the whole domain because, even if $f_u$ were not continuous on the boundary of its domain, being concave, $f_u$ is continuous in $\Intr$, and its restriction to $\Intr$ can be continuously extended to the whole domain $\sgens{k}_+$ \cite[Theorem~5.1.5]{Lemmens2012}. Such an extension of ${f_u}\restr{\Intr}$  to the domain $\sgens{k}_+$ would preserve positivity because the range of ${f_u}\restr{\Intr}$ is bounded away from zero as shown in (\ref{eq.bounded_away}).

Below we list some advantages of using the equivalence relation in (\ref{equivalent.constraint}) and the result in (\ref{eq.constraint}) to replace the inequality constraints in problem (\ref{problem.opt}) with the fixed point set of the PC mapping in (\ref{eq.int_map_pc}):
\begin{remark}
\label{remark.advantages}
\begin{itemize}
	\item[(i)] We immediately deduce that problem (\ref{problem.opt}) has at most one solution because the set $\mathrm{Fix}(f)$ is a singleton if not empty (Fact~\ref{sr1});
	\item[(ii)] Problem (\ref{problem.opt}) has a solution if and only if $\rho(f_\infty)<1$ (Fact~\ref{sr1}), which is a condition that can be easily verified with the simple numerical approaches described in \cite[Remarks~2 and 3]{Cavalcante2019_IEEE_TSP};
	\item[(iii)] If $x^\star$ is the solution to problem (\ref{problem.opt}), then $x^\star$ is the limit of the sequence $(x_n)_{n\in\mathbb{N}}$ produced via the fixed point iteration $x_{n+1}=f(x_n)$, $x_1\in\sgens{k}_+$ (Fact~\ref{yates} in the Appendix).
	\item[(iv)] The convergence of the fixed point iteration of $f$ is geometric (Proposition~\ref{p1}) if the fixed point exists.
\end{itemize}
 \end{remark}

In this particular application, the theory described in the previous sections goes beyond results available in the wireless literature w.r.t. the following aspects:

\begin{remark} \label{remark.beyond}\begin{enumerate}
		\item[(i)] Assuming that the set $\mathcal{V}_u$ is a singleton for every $u\in\mathcal{U}$, then the studies in \cite{huang1998rate,fey2012} show \emph{sufficient} conditions for the mapping $f$ in (\ref{eq.int_map_pc}) to have a fixed point. In contrast, Remark~\ref{remark.advantages}(ii) describes a simple and more general \emph{necessary and sufficient} condition for establishing the existence of the fixed point. 
		
		\item[(ii)] More importantly, if the set $\mathcal{V}_u$ contains uncountably many beamformers for some user $u\in\mathcal{U}$ (e.g., if base stations are equipped with fully digital beamformers, in which case $\mathcal{V}_u=\{v \in \mathbb{C}^L~|~{v}^H {v}=1\}$ for every $u\in\mathcal{U}$), then the mapping $f$ in (\refeq{eq.int_map_pc}) cannot be expressed as the coordinate-wise minimum of a \emph{finite} number of affine functions in general, so the theory developed in \cite{huang1998rate} cannot be used to establish geometric convergence of the fixed point iteration of $f$. Furthermore, the study in \cite{fey2012} does not describe any simple procedure to verify whether $f$ is a  $c$-contraction in the sense in \cite[Definition~2]{fey2012}, which would give us an alternative proof of geometric convergence of the fixed point iteration of $f$. In contrast to the results in \cite{fey2012, huang1998rate}, we impose no additional structure on $f$ other than positivity and concavity to establish geometric convergence of the fixed point iteration (Proposition~\ref{p1}). These properties can be immediately verified for the particular mapping in (\ref{eq.int_map_pc}) as discussed above.
		
		\item[(iii)]  Sect.~\ref{sect.specrad} shows simple techniques to obtain information about the contraction factor of the fixed point iteration of $f$ (and hence the convergence speed of the iteration) in normed vector spaces and in Thompson's metric space. Mathematical tools of this type are not available in the wireless literature to the best of our knowledge.
		\end{enumerate}
\end{remark}

We finish this section by showing that the theory developed in this study can also explain the convergence speed of extensions of the power control algorithm in Remark \ref{remark.advantages}(iii).

\begin{remark} A possible drawback of the fixed point iteration in Remark \ref{remark.advantages}(iii), which is a consequence of the optimization problem posed in (\refeq{problem.opt}), is that no maximum power constraints are imposed. To prevent users from transmitting with excessive  power, we can replace the mapping $f$ in the fixed point iteration with
\begin{align*}
	\bar{f}:\sgens{k}_+\to\Intr:x\mapsto [\min\{f_1(x),\bar{p}\},\ldots,\min\{f_u(x),\bar{p}\}],
\end{align*}
where $\bar{p}>0$ is the maximum allowed power. Note that $\bar{f}$ is also concave because each coordinate function is the minimum of concave functions. Well-known power control algorithms emerge as particular cases of the fixed point iteration of $\bar{f}$. For example, by fixing the sets $(\mathcal{V}_u)_{u\in\mathcal{U}}$ to be singletons and the random channel vectors to be almost surely constant, we obtain the power control algorithm discussed in \cite[Eq.~(7)]{sun2014} as a special case of the fixed point iteration of the resulting mapping $\bar{f}$. Whether such simplifications are used or not, we always have $\rho(\bar{f}_\infty)=0$ because each coordinate of $\bar{f}$ is upper bounded by $\bar{p}$, and thus  $\forall x\in\sgens{k}_+\quad \bar{f}_\infty(x)=0$. As a result, Fact \ref{sr1} shows that  $\mathrm{Fix}(\bar{f})\neq\emptyset$, and the fixed point iteration of $\bar{f}$ converges geometrically to the unique vector in $\mathrm{Fix}(\bar{f})$ as a consequence of Proposition~\ref{p1}.
\end{remark}



\subsection{Load estimation in orthogonal frequency division multiple access (OFDMA) networks}
\label{sect.ofdma}

As a second example of an application of the fixed point iteration of positive concave mappings, we describe the problem of load estimation in wireless networks  \cite{Majewski2010,siomina12,feh2013,ho2014data,renato14SPM,Cavalcante2016,Cavalcante2019_IEEE_TSP,cavalcante2016low,cavalcante2018spectral}. We assume that a network uses the OFDMA technology, and we consider the same downlink scenario (i.e., data transmission is from base stations to users) used in the simulations in \cite{cavalcante2018spectral}. To pose the problem, we use the following definitions:
\begin{itemize}
\item $\mathcal{B}=\{1,\ldots,k\}$ - set of base stations;
\item $\mathcal{U}$ - set of users;	
\item $\mathcal{U}_b\neq\emptyset$ - set of users connected to base station $b$;
\item $B$ - bandwidth of each resource block;
\item $\forall u\in\mathcal{U}\quad d[u]$ - traffic requested by the $u$th user;
\item $g[u,b]>0$ - pathloss between base station $b$ and user $u$;
\item $R$ - number of resource blocks in the OFDMA system;
\item $\forall b\in\mathcal{B}\quad p[b]$ - transmit power per resource block of base station $b$;
\item $\sigma^2$ - noise power per resource block.
\item $x=[x[1],\ldots, x[k]]$ - load vector. More precisely, the $b$th coordinate $x[b]$ is the fraction of the total number of resource blocks that base station $b$ requires to serve the traffic demand of its users in the set $\mathcal{U}_b$.
\end{itemize}

For a given load $x$, the modified Shannon capacity  $r_{u,b}:\sgens{k}_+\to\sgens{k}_+$, which has been independently proposed in \cite{siomina12,Majewski2010} and defined below, is often used as a proxy of the rate per resource block that a user $u\in\mathcal{U}$ can achieve if connected to base station $b\in\mathcal{B}$ \cite{Majewski2010,siomina12,feh2013,ho2014data,renato14SPM,Cavalcante2016,Cavalcante2019_IEEE_TSP,cavalcante2016low,cavalcante2018spectral}:
\begin{align*}
	\begin{array}{rl}
	r_{u,b}:\sgens{k}_+\to&\sgens{k}_+ \\ 
	x\mapsto & B\log_2\left(1+\dfrac{p[b] g[u,b]}{\sum\limits_{j\in\mathcal{B}\backslash\{b\}} x[j] p[j] g[u,j]+\sigma^2}\right).
\end{array}	
\end{align*}

The above model assumes no intracell interference; i.e., users connected to a given base station use distinct resource blocks. However, there is intercell interference, which corresponds to the term $\sum\limits_{j\in\mathcal{B}\backslash\{b\}} x[j] p[j] g[u,j]$ in the definition of $r_{u,b}$. As a result, if the load vector $x$ is given, and user $u$ is connected to base station $b$, then this base station needs to reserve $d[u]/r_{u,b}(x)$ resource blocks to satisfy the traffic demand of user $u$. By definition, the $b$th coordinate $x[b]$ of the load vector $x$ is the total number of resource blocks, normalized by number of resource blocks $R$ in the system, that base station $b$ reserves for all its users; i.e.,
\begin{align*}
	\forall b\in\mathcal{B}\quad x[b]= \dfrac{1}{R}\sum_{u\in\mathcal{U}_b}\dfrac{d[u]}{r_{u,b}(x)}.
\end{align*}
The above equations show that the load vector $x$ can be equivalently expressed as the fixed point of the following mapping:
 \begin{align} \label{eq.load_mapping}
	f:\real_+^k\to\Intr:{x}\mapsto [f_1({x}),\cdots,f_k({x})],
\end{align} 
where 
\begin{align*}
\forall b\in\mathcal{B} \quad	\begin{array}{rl}
		f_b:\sgens{k}_+\to&\Intr \\
		x\mapsto& \dfrac{1}{R}\sum_{u\in\mathcal{U}_b}\dfrac{d[u]}{r_{u,b}(x)}.	
	\end{array}
\end{align*}

Therefore, to estimate the load $x$ of a network, we need to solve the following fixed point problem:
\begin{align}
	\label{eq.load_est}
	\text{Find }x\in\sgens{k}_+ \text{ such that } x\in\mathrm{Fix}(f).
\end{align}

Standard arguments in convex analysis show that $f$ is a continuous positive concave mapping \cite{siomina12,Majewski2010,renato14SPM}. Therefore, all the theory in the previous sections can be applied to analyze and solve problem (\ref{eq.load_est}). In particular, we deduce the following properties:
\begin{enumerate}
	\item[(i)] Problem~(\ref{eq.load_est}) has at most one solution (Fact~\ref{sr1});
	\item[(ii)] The asymptotic mapping associated with $f$ is given by \cite{Cavalcante2019_IEEE_TSP}
	$$f_\infty:\sgens{k}_+\to\sgens{k}_+:{x}\mapsto \mathrm{diag}({p})^{-1}{M}\mathrm{diag}({p}){x},$$
	where $\mathrm{diag}({p})$ is a diagonal matrix with diagonal elements given by $\forall b\in\mathcal{B}\quad[\mathrm{diag}({p})]_{bb}=p[b]$, and $M\in\sgens{k\times k}_+$ is the matrix with its $i$th row and $b$th column given by:
	\begin{align*}
		[{M}]_{i,b} = \begin{cases}
			0,&\text{if }i=b\\
			\sum\limits_{u \in \mathcal{U}_i} \dfrac{\mathrm{ln}(2)  d[u] g[u,b]}{{RB g[u,i]}}&\text{otherwise.}
		\end{cases}
	\end{align*}
	Note that the asymptotic mapping $f_\infty$ is linear, so $\rho(f_\infty)=\rho(\mathrm{diag}({p})^{-1}{M}\mathrm{diag}({p}))=\rho(M)$ is simply the spectral radius $\rho(M)=\rho(\mathrm{diag}({p})^{-1}{M}\mathrm{diag}({p}))$ of the matrix ${M}$, where the spectral radius of matrices should be understood in the conventional sense in linear algebra. It now follows from Fact~\refeq{sr1} that
	\begin{align*}
		\mathrm{Fix}(f)\neq\emptyset \iff \rho(M)<1.
	\end{align*}
	This criterion for checking the existence of a feasible load vector has been first obtained in \cite{ho2014data}, and it can be seen as a particular application of Fact~\refeq{sr1}, as discussed in \cite{Cavalcante2019_IEEE_TSP}.
	\item[(iii)] If $\mathrm{Fix}(f)\neq\emptyset$, then the unique fixed point of $f$ is the limit of the sequence generated via the fixed point iteration of $f$ (Fact~\ref{yates} in the Appendix). 
    \item[(iv)] The fixed point iteration of $f$ converges geometrically (Proposition~\ref{p1}), and the convergence factor $c$ of the iteration is lower bounded by the spectral radius of the matrix $M$ defined above (Proposition~\ref{prop.linearc}). 
\end{enumerate}

We emphasize that points (i)-(iii) are known, but point (iv) is a new result that does not seem easy to derive from existing theory for reasons similar to those given in Remark~\ref{remark.beyond}. More precisely, the mapping in (\ref{eq.load_mapping}) cannot be written as the coordinate-wise minimum of a finite number of affine mappings, as required in \cite{huang1998rate}. Furthermore, it is not clear whether that mapping satisfies the conditions in \cite[Definition~2]{fey2012}. Even if this last condition could be established, the study in \cite{fey2012} does not discuss any simple computational means, such as those described in Sect.~\ref{sect.specrad}, to determine bounds on the contraction factor.

Before we proceed with numerical examples to illustrate the main findings in this study, we recall one potential caveat of the load model discussed here:

\begin{remark}
	By definition, the coordinates of the load vector $x$ cannot exceed the value one, otherwise the corresponding base stations would be using more resource blocks than available in the system, which is physically impossible. Nevertheless, if any coordinate of the solution to problem (\ref{eq.load_est}) exceeds the value one, we obtain useful information about the unserved traffic demand. We refer the readers to \cite{siomina12} for more information on this point.
\end{remark}

\section{Simulations}
\label{simulations}

For brevity, we show numerical examples only for the load estimation problem described in Sect.~\ref{sect.ofdma}. In the simulation, the pathloss $g[u,b]$ between base station $b\in\mathcal{B}$ and user $u\in\mathcal{U}$ is computed using the standard Hata model for urban scenarios with the height of the antennas of the users and the base stations set to, respectively, $1.5$m and $30$m. The number of base stations is $k=25$, and the base stations are uniformly distributed in a square of dimensions $2\text{km}\times2\text{km}$. A total of $400$ users are distributed uniformly at random in a square of dimensions $2.5\text{km}\times 2.5\text{km}$, which is concentric to the square where base stations are placed. Users are assigned to base stations with the lowest pathloss. Other main parameters of the simulations are listed in Table~\ref{table.parameters}.

\begin{table}
\begin{center}
\begin{tabular}{l l}
	\hline
	Parameter & Value \\
	\hline 
    $R$ &  25 \\
	$B$ &   $2\cdot 10^5$~Hz \\
    $\forall u\in\mathcal{U}\quad d[u]$ & $10^6$ bits/second \\
	$\forall b\in\mathcal{B}\quad p[b]$ &  1.6~W\\
    $\sigma^2$ & $6.2\cdot 10^{-18}$~W\\
	\hline
\end{tabular}
\caption{Main parameters used in the simulation.}
\label{table.parameters}
\end{center}
\end{table}

\begin{figure}
	\begin{center}
		\includegraphics[width=\columnwidth]{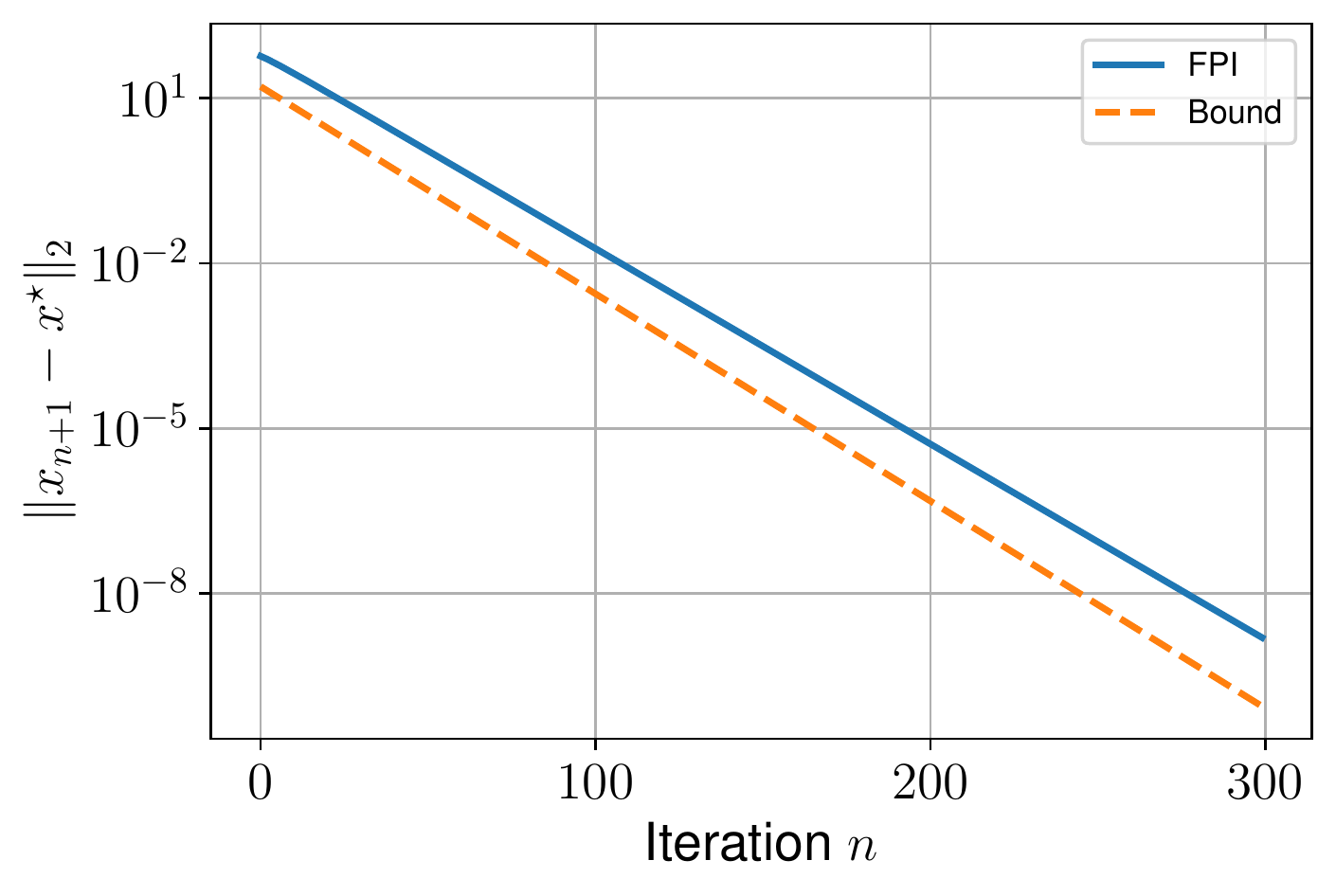}
		\caption{Convergence of the fixed point iteration as a function of the number of iterations.}
		\label{fig.convergence}
	\end{center}
\end{figure}

\begin{figure}
	\begin{center}
		\includegraphics[width=\columnwidth]{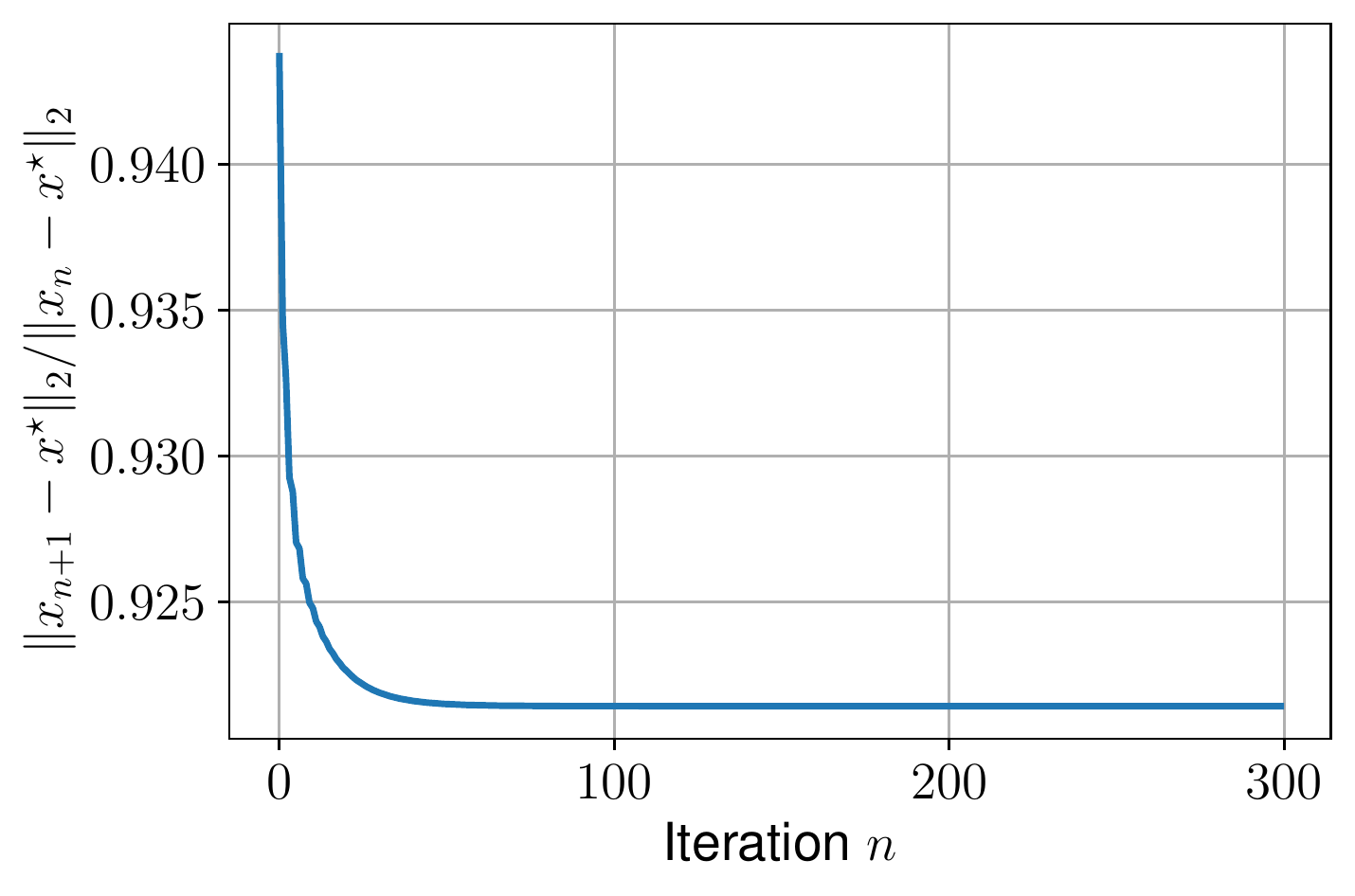}
		\caption{Ratio ${\|x_{n+1}-x^\star\|_2}/{\|x_n-x^\star\|_2}$ as a function of the iteration index~$n$.}
		\label{fig.rate}
	\end{center}
\end{figure}

For the simulated scenario, we have $\rho(M)\approx 0.92,$ so the fixed point exists, and it has been computed by constructing the sequence $(x_n)_{n\in\mathbb{N}}$ via the fixed point iteration $x_{n+1}=f(x_n)$, where $x_1=f(0)$. In Fig.~\ref{fig.convergence} we show the convergence of the fixed point iterations and also the lower bound derived in \cite[Proposition~5(iii)]{cavalcante2018spectral}. In Fig.~\ref{fig.rate} we show the ratio ${\|x_{n+1}-x^\star\|_2}/{\|x_n-x^\star\|_2}$ as a mapping of the number of iterations. In these figures, $\|\cdot\|_2$ denotes the standard Euclidean norm. 

In light of Proposition~\ref{prop.linearc}, $\rho(M)$ is a lower bound for the local contraction factor around the fixed point $x^\star$ of $f$, and this property is indeed consistent with the result in Figs.~\ref{fig.convergence} and \ref{fig.rate}. More specifically, note that the bound has a slightly greater slope than the error $\|x_n-x^\star\|_2$ for the iterations shown in Fig.~\ref{fig.convergence}. This property is expected because the bound as a mapping of the iteration index $n$ takes the form $y[n] = C_1 \rho(M)^n$ for some $C_1>0$, whereas the error takes the form $\|x_{n+1}-x^\star\|_2\le C_2 c^n$ for some $C_2>0$ and some $c\in [\rho(M),1)$. Furthermore, the ratio in Fig.~\ref{fig.rate} is converging to a value $\mu$ within the interval $(0,~1)$, which provides us with an indication of the type of convergence described in Proposition~\ref{prop.linearc} with $\mu\ge\rho(M)$ \cite{luke2020necessary}. Note, however, that the results in that proposition are valid for any sequence produced via the fixed point iteration, and this simulation shows the convergence of a sequence produced by starting the iterations with a particular vector.

\section{Summary and final remarks}
Previous studies in the wireless literature \cite{Cavalcante2016,cavalcante2018spectral,fey2012,huang1998rate}\cite[Ch.~5.3]{martin11} have devoted a great deal of attention to characterize the convergence speed of fixed point iterations of mappings that are common in power control and load estimation problems. In many of the these applications, the mappings are in fact positive and concave. However, as discussed in Sect.~\ref{applications}, previous results proving geometric convergence of the fixed point iteration of positive concave mappings require additional assumptions that are often too strong to address some problems of practical interest, or that are too difficult to verify in practice, or both. In this study we have proved that these additional assumptions are not required: we only need the existence of a fixed point for the fixed point iteration of positive concave mappings to converge geometrically. We have also shown that even by replacing the assumption of positivity with nonnegativity we lose the guarantees of geometric convergence, so the assumption of positivity is crucial. Finally, we have also derived results that provide us with clear indication of the local contraction factor of the mappings, which is useful information to assess whether many fixed point iterations are required to obtain a good estimate of the fixed point.

\appendices
\section{Definitions and Known Results} \label{kru}
\begin{definition}(\cite[p.~45]{Lemmens2012}) \label{nc}
Let $\|\cdot\|$ be any norm on $\sgens{k}.$  We call $\delta>0$ the normality constant of the cone $\sgens{k}_+$ for the norm $\|\cdot\|$ if
\begin{equation}
\delta\df\inf\{q\in\mathbb{R}_+~|~\forall(x,y)\in C\quad\|x\|\leq q\,\|y\|\},
\end{equation}
where
$C\df\{(x,y)\in\sgens{k}_+\times\sgens{k}_+~|~ x\leq y\}.$
\end{definition}

\begin{fact}[$\mbox{\cite[Lemma 2.5.1.]{Lemmens2012}}$] \label{L251} 
Let $K$ be a cone in a normed space $(V,\|\cdot\|),$\footnote{For example, $V=\sgens{k}$ with $K=\sgens{k}_+.$} and let $\delta>0$ be the normality constant of~$K.$ For each $x,y\in K$ with $\|x\|\leq b$ and $\|y\|\leq b$, we have that
\begin{equation*}
\|x-y\|\leq b(1+2\delta)(e^{d_T(x,y)}-1).
\end{equation*}  
\end{fact}

\begin{fact}[$\mbox{\cite[Theorem 2]{yates95}}$] \label{yates} 
Let $f\colon\sgens{k}_+\to\Intr$ be a continuous SI mapping satisfying the feasibility condition $f(x)\leq x$ for a vector $x\in\sgens{k}_+.$ Then, for any $x\in\Intr$, the fixed point iteration of $f$ produces a sequence $(x_n)_{n\in\mathbb{N}}$ converging in norm to the unique vector $x^\star\in\mathrm{Fix}(f)$, i.e., $\lim_{n\to\infty}\|x_n-x^\star\|=0$ for an arbitrary norm $\|\cdot\|.$ 
\end{fact}

\bibliographystyle{IEEEtran}
\bibliography{IEEEabrv,references}

\begin{thebibliography}{10}
\providecommand{\url}[1]{#1}
\csname url@samestyle\endcsname
\providecommand{\newblock}{\relax}
\providecommand{\bibinfo}[2]{#2}
\providecommand{\BIBentrySTDinterwordspacing}{\spaceskip=0pt\relax}
\providecommand{\BIBentryALTinterwordstretchfactor}{4}
\providecommand{\BIBentryALTinterwordspacing}{\spaceskip=\fontdimen2\font plus
\BIBentryALTinterwordstretchfactor\fontdimen3\font minus
  \fontdimen4\font\relax}
\providecommand{\BIBforeignlanguage}[2]{{%
\expandafter\ifx\csname l@#1\endcsname\relax
\typeout{** WARNING: IEEEtran.bst: No hyphenation pattern has been}%
\typeout{** loaded for the language `#1'. Using the pattern for}%
\typeout{** the default language instead.}%
\else
\language=\csname l@#1\endcsname
\fi
#2}}
\providecommand{\BIBdecl}{\relax}
\BIBdecl

\bibitem{woodford2020}
M.~Woodford, ``Effective demand failures and the limits of monetary
  stabilization policy,'' National Bureau of Economic Research, Tech. Rep.,
  2020.

\bibitem{yates95}
R.~D. Yates, ``A framework for uplink power control in cellular radio
  systems,'' \emph{{IEEE} J. Sel. Areas Commun.}, vol.~13, no.~7, pp. pp.
  1341--1348, Sept. 1995.

\bibitem{Cavalcante2016}
R.~L.~G. Cavalcante, Y.~Shen, and S.~Stańczak, ``Elementary properties of
  positive concave mappings with applications to network planning and
  optimization,'' \emph{IEEE Transactions on Signal Processing}, vol.~64,
  no.~7, pp. 1774--1783, 2016.

\bibitem{Cavalcante2019_IEEE_TSP}
R.~L.~G. Cavalcante, Q.~Liao, and S.~Stańczak, ``Connections between spectral
  properties of asymptotic mappings and solutions to wireless network
  problems,'' \emph{IEEE Transactions on Signal Processing}, vol.~67, no.~10,
  pp. 2747--2760, 2019.

\bibitem{martin11}
M.~Schubert and H.~Boche, \emph{Interference Calculus - A General Framework for
  Interference Management and Network Utility Optimization}.\hskip 1em plus
  0.5em minus 0.4em\relax Berlin: Springer, 2011.

\bibitem{slawomir09}
S.~Sta\'nczak, M.~Wiczanowski, and H.~Boche, \emph{Fundamentals of Resource
  Allocation in Wireless Networks}, 2nd~ed., ser. Foundations in Signal
  Processing, Communications and Networking, W.~Utschick, H.~Boche, and
  R.~Mathar, Eds.\hskip 1em plus 0.5em minus 0.4em\relax Berlin Heidelberg:
  Springer, 2009.

\bibitem{you2020note}
L.~You and D.~Yuan, ``A note on decoding order in user grouping and power
  optimization for multi-cell {NOMA} with load coupling,'' \emph{IEEE
  Transactions on Wireless Communications}, 2020.

\bibitem{shindoh2020}
S.~Shindoh, ``Some properties of {SINR} regions for standard interference
  mappings,'' \emph{SICE Journal of Control, Measurement, and System
  Integration}, vol.~13, no.~3, pp. 50--56, 2020.

\bibitem{cavalcante2016low}
R.~L.~G. Cavalcante, S.~Sta{\'n}czak, J.~Zhang, and H.~Zhuang, ``Low complexity
  iterative algorithms for power estimation in ultra-dense load coupled
  networks,'' \emph{IEEE Transactions on Signal Processing}, vol.~64, no.~22,
  pp. 6058--6070, 2016.

\bibitem{ho2014data}
C.~K. Ho, D.~Yuan, and S.~Sun, ``Data offloading in load coupled networks: A
  utility maximization framework,'' \emph{IEEE Transactions on Wireless
  Communications}, vol.~13, no.~4, pp. 1921--1931, 2014.

\bibitem{piotrowski2021fixed}
T.~Piotrowski and R.~L.~G. Cavalcante, ``Fixed points of nonnegative neural
  networks,'' \emph{submitted, also available as arXiv preprint
  arXiv:2106.16239}, 2021.

\bibitem{Lemmens2012}
B.~Lemmens and R.~Nussbaum, \emph{Nonlinear Perron-Frobenius Theory}.\hskip 1em
  plus 0.5em minus 0.4em\relax Cambridge Univ. Press, 2012.

\bibitem{Oshime1992}
Y.~Oshime, ``{Perron-Frobenius} problem for weakly sublinear maps in a
  {Euclidean} positive orthant,'' \emph{Japan Journal of Industrial and Applied
  Mathematics}, vol.~9, no. 313, 1992.

\bibitem{krause2015positive}
U.~Krause, \emph{Positive dynamical systems in discrete time: theory, models,
  and applications}.\hskip 1em plus 0.5em minus 0.4em\relax Walter de Gruyter
  GmbH \& Co KG, 2015, vol.~62.

\bibitem{Majewski2010}
K.~Majewski and M.~Koonert, ``Conservative cell load approximation for radio
  networks with {Shannon} channels and its application to {LTE} network
  planning,'' in \emph{Telecommunications (AICT), 2010 Sixth Advanced
  International Conference on}, May 2010, pp. 219 --225.

\bibitem{siomina12}
I.~Siomina and D.~Yuan, ``Analysis of cell load coupling for {LTE} network
  planning and optimization,'' \emph{{IEEE} Trans. Wireless Commun.}, vol.~11,
  no.~6, pp. 2287--2297, June 2012.

\bibitem{feh2013}
A.~Fehske, H.~Klessig, J.~Voigt, and G.~Fettweis, ``Concurrent load-aware
  adjustment of user association and antenna tilts in self-organizing radio
  networks,'' \emph{{IEEE} Trans. Veh. Technol.}, no.~5, June 2013.

\bibitem{shindoh2019}
S.~Shindoh, ``The structures of {SINR} regions for standard interference
  mappings,'' in \emph{2019 58th Annual Conference of the Society of Instrument
  and Control Engineers of Japan (SICE)}.\hskip 1em plus 0.5em minus
  0.4em\relax IEEE, 2019, pp. 1280--1285.

\bibitem{cavalcante2018spectral}
R.~L.~G. Cavalcante and S.~Sta{\'n}czak, ``Spectral radii of asymptotic
  mappings and the convergence speed of the standard fixed point algorithm,''
  in \emph{2018 IEEE International Conference on Acoustics, Speech and Signal
  Processing (ICASSP)}.\hskip 1em plus 0.5em minus 0.4em\relax IEEE, 2018, pp.
  4509--4513.

\bibitem{fey2012}
H.~R. Feyzmahdavian, M.~Johansson, and T.~Charalambous, ``Contractive
  interference functions and rates of convergence of distributed power control
  laws,'' \emph{IEEE Transactions on Wireless Communications}, vol.~11, no.~12,
  pp. 4494--4502, 2012.

\bibitem{huang1998rate}
C.-Y. Huang and R.~D. Yates, ``Rate of convergence for minimum power assignment
  algorithms in cellular radio systems,'' \emph{Wireless Networks}, vol.~4,
  no.~4, pp. 223--231, 1998.

\bibitem{boche2008superlinearly}
H.~Boche and M.~Schubert, ``A superlinearly and globally convergent algorithm
  for power control and resource allocation with general interference
  functions,'' \emph{IEEE/ACM Transactions on Networking}, vol.~16, no.~2, pp.
  383--395, 2008.

\bibitem{hanly1995algorithm}
S.~V. Hanly, ``An algorithm for combined cell-site selection and power control
  to maximize cellular spread spectrum capacity,'' \emph{IEEE Journal on
  selected areas in communications}, vol.~13, no.~7, pp. 1332--1340, 1995.

\bibitem{luke2020necessary}
D.~R. Luke, M.~Teboulle, and N.~H. Thao, ``Necessary conditions for linear
  convergence of iterated expansive, set-valued mappings,'' \emph{Mathematical
  Programming}, vol. 180, no.~1, pp. 1--31, 2020.

\bibitem{Lauster2021}
F.~Lauster and D.~R. Luke, ``{Convergence of proximal splitting algorithms in
  $\operatorname{CAT}(\kappa)$ spaces and beyond},'' \emph{Fixed Point Theory
  and Algorithms for Sciences and Engineering}, vol. 2021, no.~13, 2021.

\bibitem{boche2008}
H.~Boche and M.~Schubert, ``Concave and convex interference functions --
  general characterizations and applications,'' \emph{IEEE Transactions on
  Signal Processing}, vol.~56, no.~10, pp. 4951--4965, 2008.

\bibitem{Ortega1970}
J.~M. Ortega and W.~C. Rheinboldt, \emph{Iterative Solution of Nonlinear
  Equations in Several Variables}.\hskip 1em plus 0.5em minus 0.4em\relax New
  York: Academic Press, 1970.

\bibitem{gau04}
S.~Gaubert and J.~Gunawardena, ``The {Perron-Frobenius} theorem for
  homogeneous, monotone functions,'' \emph{Transactions of the American
  Mathematical Society}, vol. 356, no.~12, pp. 4931--4950, 2004.

\bibitem{Bauschke2017}
H.~H. Bauschke and P.~L. Combettes, \emph{Convex Analysis and Monotone Operator
  Theory in Hilbert Spaces}.\hskip 1em plus 0.5em minus 0.4em\relax New York:
  Springer, 2017.

\bibitem{Leung2004}
K.~K. Leung, C.~W. Sung, W.~S. Wong, and T.-M. Lok, ``Convergence theorem for a
  general class of power-control algorithms,'' \emph{{IEEE} Trans. Commun.},
  vol.~52, no.~9, pp. 1566--1574, 2004.

\bibitem{Cavalcante2019}
R.~L.~G. {Cavalcante} and S.~{Stańczak}, ``Weakly standard interference
  mappings: Existence of fixed points and applications to power control in
  wireless networks,'' in \emph{ICASSP 2019 - 2019 IEEE International
  Conference on Acoustics, Speech and Signal Processing (ICASSP)}, 2019, pp.
  4824--4828.

\bibitem{Potter1977}
A.~J.~B. Potter, ``{Applications of Hilbert's projective metric to certain
  classes of non-homogeneous operators},'' \emph{The Quarterly Journal of
  Mathematics}, vol.~28, no.~1, pp. 93--99, 1977.

\bibitem{Chengbo2005}
Z.~Chengbo and G.~Chunmei, ``{On \textalpha-convex operators},'' \emph{Journal
  of Mathematical Analysis and Applications}, vol. 316, pp. 556--565, 2005.

\bibitem{Liang2006}
Z.~D. Liang, W.~X. Wang, and S.~J. Li, ``{On concave operators},'' \emph{Acta
  Mathematica Sinica, English Series}, vol.~22, no.~2, pp. 577--582, 2006.

\bibitem{BabyRudin}
W.~Rudin, \emph{Principles of Mathematical Analysis}.\hskip 1em plus 0.5em
  minus 0.4em\relax New York: McGraw-Hill, 1976.

\bibitem{Nussbaum1986}
R.~D. Nussbaum, ``Convexity and log convexity for the spectral radius,''
  \emph{Linear Algebra and its Applications}, vol.~73, pp. 59--122, 1986.

\bibitem{sun2014}
R.~Sun and Z.-Q. Luo, ``Globally optimal joint uplink base station association
  and power control for max-min fairness,'' in \emph{IEEE International
  Conference on Acoustics, Speech and Signal Processing (ICASSP)}, 2014, pp.
  454--458.

\bibitem{renato14SPM}
R.~L.~G. Cavalcante, S.~Sta\'nczak, M.~Schubert, A.~Eisenbl\"ater, and
  U.~T\"urke, ``Toward energy-efficient {5G} wireless communication
  technologies,'' \emph{{IEEE} Signal Process. Mag.}, vol.~31, no.~6, pp.
  24--34, Nov. 2014.

\end{thebibliography}

\end{document}